\documentclass{article}

\usepackage{arxiv}

\usepackage[utf8]{inputenc} 
\usepackage[T1]{fontenc}    
\usepackage{hyperref}       
\usepackage{url}            
\usepackage{booktabs}       
\usepackage{amsfonts}       
\usepackage{nicefrac}       
\usepackage{microtype}      
\usepackage{lipsum}
\usepackage{algorithm}
\usepackage{algorithmic}
\usepackage{amsmath}
\usepackage{amssymb}
\usepackage{epstopdf}
\usepackage{enumitem}
\usepackage{float}
\usepackage{graphicx}
\usepackage{lineno,hyperref}
\usepackage{lscape}
\usepackage{morefloats}
\usepackage{multicol}
\usepackage{multirow}
\usepackage{subfigure}
\usepackage{ulem}
\usepackage{longtable}
\usepackage{multirow}
\usepackage{array}
\usepackage[table]{xcolor}
\usepackage{threeparttable}

\usepackage{color}
\newcommand{\change}[1]{\textcolor{black}{#1}}  
\newcommand{\blue}[1]{\textcolor{black}{#1}}

\def\bfa{{\mathbf{a}}}
\def\bfb{{\mathbf{b}}}

\def\bfg{{\mathbf{g}}}

\def\bfu{{\mathbf{u}}}
\def\bfv{{\mathbf{v}}}
\def\bfx{{\mathbf{x}}}

\def\bfz{{\mathbf{z}}}
\def\bfdelta{{\boldsymbol{\delta}}}

\def\bfphi{{\boldsymbol{\phi}}}

\def\bfxi{{\boldsymbol{\xi}}}

\def\du{\textrm{d}\bfu}

\def\dxi{\textrm{d}\bfxi}

\def\calG{\mathcal{G}}

\def\calU{\mathcal{U}}

\def\calJ{\mathcal{J}}
\def\calL{\mathcal{L}}

\title{Parametric Model Embedding}

\author{
  Andrea Serani\thanks{www.inm.cnr.it} \\
  CNR-INM, National Research Council\\
  Institute of Marine Engineering\\
  Rome, Italy  \\
  \texttt{andrea.serani@cnr.it} \\
   \And
 Matteo Diez \\
  CNR-INM, National Research Council\\
  Institute of Marine Engineering\\
  Rome, Italy  \\
  \texttt{matteo.diez@cnr.it} \\
}

\begin{document}
\maketitle

\begin{abstract}
Methodologies for reducing the design-space dimensionality in shape optimization have been recently developed based on unsupervised machine learning methods. These methods provide reduced dimensionality representations of the design space, capable of maintaining a certain degree of the original design variability. Nevertheless, they usually do not allow to use \change{directly} the original parameterization \change{method}, representing a limitation to their widespread application in the industrial field, where the design parameters often pertain to well-established parametric models, e.g. CAD (computer aided design) models. This work presents how to embed the parametric-model original parameters in a reduced-dimensionality representation \change{of the design space}. The method, which takes advantage from the definition of a newly-introduced generalized feature space, is demonstrated, \blue{as a proof of concept,} for the reparameterization of \blue{2D Bezier curves and 3D} free-form deformation design space\blue{s} and the consequent solution of simulation-driven design optimization problem\blue{s} of \blue{a subsonic airfoil and} a naval destroyer in calm water, \blue{respectively}.
\end{abstract}


\change{
\section*{Acronyms and nomenclature}
\begin{table}[!h]
\centering
\begin{tabular}{lcllcl}
\toprule
CoD &= & curse of dimensionality                & $\bfg$ &=     & original (parent) shape\\
FFD &= & free-form deformation                  & $\bfdelta$ &= & shape modification vector\\
KLE &= & Karhune-Lo{\`e}ve expansion            & $\bfu$ &=     & original design-variables vector\\
PCA &= & principal component analysis           & $\bfx$ &=     & reduced design-variables vector\\
PME &= & parametric model embedding             & $\bfz$ &=     & geometrical eigenvector component\\
POD &= & proper orthogonal decomposition        & $\bfv$ &=     &  parameterization eigenvector component\\ 
SDDO &= & simulation-driven design optimization & $M$ &= & number of original design variables\\
$\bfxi$ &= & shape curvilinear coordinates      & $N$ &= & number of reduced design variables\\
\bottomrule
\end{tabular}
\end{table}
}

\section{Introduction}
The need for increasingly performing functional-surface designs is constantly growing in many engineering fields, requiring increasingly accurate analysis and innovative solutions. The latter can be achieved via the simulation-driven design optimization (SDDO) paradigm \cite{harries2019-OE}, which integrates shape parameterization models, numerical solvers for the physics of interest, and optimization algorithms. The demand for highly innovative designs often requires global optimization on ever-larger design spaces, with an ever-increasing number of design variables, leading unavoidably to the so-called curse of dimensionality (CoD) \cite{bellman1957-CS}, for which the performance of an optimization algorithm degrades as the dimensionality of the problem increases. This strongly motivates the need for reducing the dimensionality of the design space before the optimization \change{is actually performed}, especially in industrial design, where time resources are generally limited \cite{serani2021-EWCO}. Historically, the simplest method to reduce the dimensionality of the design space is to identify the most important variables for the design problem and discard the remaining ones by setting them to a constant value during the optimization process, i.e. a factor screening \cite{montgomery1979-IEEE}, also know as feature selection. However, this approach does not always provide the best solution, as it is not able to evaluate the importance that the fixed variables could have during the optimization process, especially when combined with other variables. \change{This can be assessed by variance-based sensitivity analysis (i.e., the Sobol indices \cite{sobol2001global}), nevertheless it works in a probabilistic framework and consequently needs a statistically significant number of design-space samples; furthermore, the number of indices increase with the power of the design-space dimensionality, so Sobol indices are generally computationally expensive.}
Hence, the need in industrial design is for such dimensionality reduction methods that can capture, in a reduced-dimensionality space, the underlying most promising directions of the original design space, preserving its relevant features and thereby enabling an efficient and effective optimization in the reduced space.
The remedy can be found in those dimensionality reduction techniques classified as unsupervised learning, feature extraction, or also representation learning \cite{bengio2013-IEEE}, capable of learning relevant hidden structures of the original design-space parameterization.



High-dimensional SDDO problems in functional-surface design may be reduced in dimensionality using off-line (or upfront) design-space dimensionality reduction methods. These methods have been developed focusing on the assessment of design-space variability and the subsequent dimensionality reduction before the optimization is performed. A method based on the Karhunen-Loève expansion (KLE, equivalent to the proper orthogonal decomposition, POD) has been formulated by Diez et al. \cite{diez2015-CMAME} for the assessment of the shape modification variability and the definition of a reduced-dimensionality global model of the shape modification vector. No objective function evaluation nor gradient is required by the method. The KLE is applied to the continuous shape modification vector, requiring the solution of 
\change{an eigenproblem of an integral operator}.
Once the equation is discretized, the problem reduces to the principal component analysis (PCA) of discrete geometrical data. Off-line methods improve the shape optimization efficiency by reparameterization and dimensionality reduction, providing the assessment of the design space and the shape parametrization before optimization and/or performance analysis are carried out. The assessment is based on the geometric variability associated to the design space, making the method computationally very efficient and attractive (no simulations are required). This method and closely related approaches, \change{also known as modal parameterization \cite{li2022-PAS}}, have been successfully applied to several SDDO problems, e.g. for airfoil \cite{poole2017-CF,cinquegrana2018-CF,yasong2018-CJA}, \change{wing \cite{allen2018wing}, rotor blade \cite{yanhui2019performance},} marine propeller \cite{gaggero2020-SOS}, hydropower turbine \cite{masood2021-RE}, hull-form \cite{tezzele2018-AMSES,dagostino2020-OE,harries2021application,liu2021-OE,ccelik2021reduced}, and \change{optics} \cite{melati2019mapping,torrijos2021design} applications, \blue{where such methodologies have demonstrated their capability in reducing the design-space dimensionality before the optimization loop, alleviating the CoD associated with SDDO problems \cite{diez2016-AIAA}. Nevertheless, if the dimensionality reduction procedure is fed only with information on the shape modification vector, these methods does not directly provide a way to return to the original design variables from the so-called latent space (i.e., the reduced dimensionality representation of the original shape parameterization)}.
\blue{Two significant criticalities ensue: (1) shape modification-based KLE/PCA obliges the user to implement a new shape modification method based on the KLE/PCA eigenfunctions; (2) moreover, depending on the bounds applied to the reduced design variables, there is not guarantee that the shape produced using KLE/PCA eigenvectors actually belongs to the original design space, thus potentially resulting in design unfeasibilities.}

A rigorous method able to provide an explicit relationship between the original and reduced design variables still needs to be found. This aspect represents a limit to the widespread use of \blue{KLE/PCA-based design-space dimensionality reduction methods} in the industrial field, where the design parameters often pertain to well-established parametric models (e.g., computer aided design, CAD, models), which are often used as a required input to the numerical solvers {(e.g., for isogeometric analysis)}.
To this end, the challenge is to extend the methodology to solve the so-called pre-image problem \cite{gaudrie2020-SMO}, where the objective is to build a direct mapping of the original variables from the reduced design space, i.e. to build a reduced-dimensionality embedding of the original design parameterization. To this aim\change{, to the authors best knowledge, only few recent works have been proposed to solve the pre-image problem, based on interpolation/regression methods: a back-transformation approach has been proposed in \cite{bergmann2018massive}, based on the interpolation between points in KLE-space for which are know the parameters in the original parametric space; Gaudrie et al. \cite{gaudrie2020-SMO} have proposed a methodology relying on a Gaussian process model. Nevertheless, interpolation/regression approaches are not able to provide any explicit mathematical relationships between the original parameterization and its  reduced-dimensionality representation and/or their application may be not straightforward.}  

The objective of this work is to \change{introduce and} discuss \change{a novel approach able} to easily embed \blue{(or contain)} the parametric-model original parameters in a reduced-dimensionality representation. The method, introduced here as parametric model embedding (PME), provides a direct way to return to the original parameterization from the reduced-dimensionality latent space, allowing therefore for the super-parametrization of the original parametric model. \change{These \textit{super parameters} exist in the reduced-dimensionality space and are not just a subset of the original variables in the parametric space, as it happens in a classical sensitivity analysis, where the most important variables are identified from a design of experiments \cite{harries2019-OE}.}
Specifically, the proposed PME method extends the formulation presented in \cite{diez2015-CMAME} using a generalized feature space that includes shape modification and design variables vectors together with a generalized inner product, aiming at resolving a prescribed design variability by properly selecting the latent dimensionality. \blue{As a proof of concept,} a demonstration of the PME method is provided for \blue{the design-space reparameterization and the consequent solution of two single-objective SDDO problems: (1) a 14-design-variables space based on Bezier curves for the drag coefficient reduction of an airfoil in subsonic condition and (2) a 22-design-variables space based on free-form deformation (FFD) \cite{sederberg1986-CGIT} for the resistance reduction of a destroyer-type vessel in calm water.}

\change{The remainder of the paper is organized as follows. Section \ref{sec:II} introduces the notation used and recall the design-space dimensionality reduction by KLE as per \cite{diez2015-CMAME}; its extension to PME method is presented and discussed in Section \ref{sec:III}; their application to the SDDO problems and the numerical
results are presented in Section \ref{sec:IV}. The conclusions and possible future works are discussed in Section \ref{sec:V} and, finally, \ref{sec:appA} includes the demonstration of the relationship between KLE and PME eigensolutions.}

\section{Design-space dimensionality reduction \change{for shape optimization}}\label{sec:II}
\begin{figure}[!b]
\centering
\includegraphics[width=0.95\textwidth]{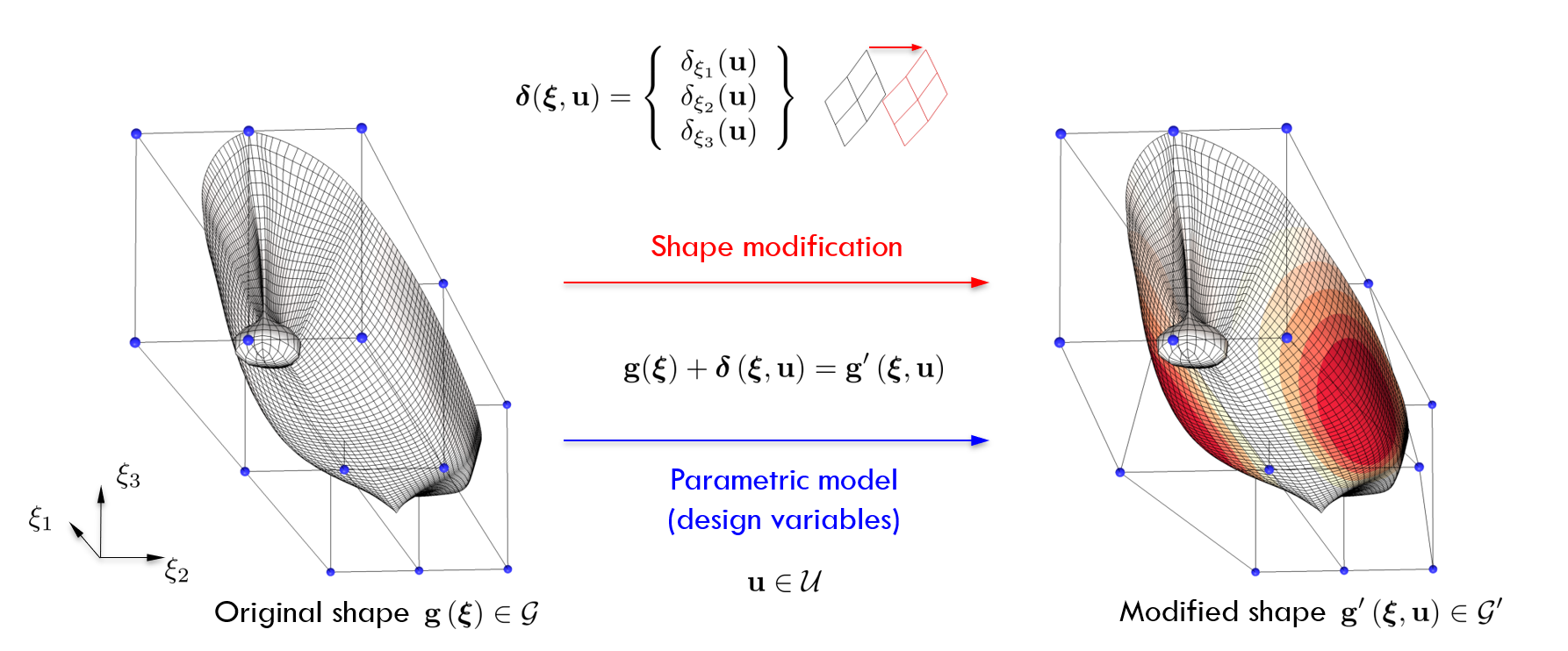}
\caption{Shape modification example and notation}\label{fig:notation}
\end{figure}

Consider a manifold $\calG$, which identifies the original/parent shape, whose coordinates in the 3D-space are represented by $\bfg(\bfxi)\in\mathbb{R}^3$; $\bfxi\in\calG$ are curvilinear coordinates defined on $\calG$. Assume that, for the purpose of shape optimization, $\bfg$ can be transformed to a deformed shape/geometry $\bfg'(\bfxi,\bfu)$ by
\begin{equation}
\bfg'(\bfxi,\bfu)=\bfg(\bfxi)+\bfdelta(\bfxi,\bfu) \qquad \forall \bfxi\in\calG
\end{equation} 
where $\bfdelta(\bfxi,\bfu)\in\mathbb{R}^3$ is the resulting shape modification vector, defined by arbitrary shape parameterization or modification method (e.g., CAD parameterization, Bezier surfaces, FFD, NURBS, etc.), and $\bfu\in\calU\subset\mathbb{R}^M$ is the design variable vector. Figure \ref{fig:notation} shows an example of the current notation.

Consider, temporarily, the \change{identification} of the optimal design, through a SDDO optimization problem
\begin{equation}
\underset{\bfu\in\calU}{\rm minimize}\,\,\, f(\bfu)
\end{equation}
as a problem affected by epistemic uncertainty, where $\bfu$ can be assumed as an uncertain/random parameter. The idea is that the optimal design exists but, before going through the optimization procedure, is unknown. Accordingly, the design variable vector $\bfu$ is assigned with a probability \change{density} function (PDF) $p(\bfu)$, representing the degree of belief in finding the optimal solution in \change{a} certain region of the design space.
\change{This probability density is arbitrary and may be chosen based on the designer's previous knowledge and past experience. If previous knowledge is not available or deemed not useful for the current optimization study, a uniform distribution may be used, which gives the same probability to all designs in the given domain to be selected as optimal decision.}

Once $p(\bfu)$ is defined, the shape modification vector $\bfdelta$ goes stochastic and can be studied as a random field, e.g., by using the KLE, equivalent to POD.

\subsection{Karhunen-Lo{\`e}ve expansion: dimensionality reduction in the continuum geometric space}\label{sec:kle}
Consider $\bfdelta(\bfxi,\bfu)$ as belonging to a Hilbert space $L_\rho^2(\calG)$, defined by the generalized inner product
\begin{equation}\label{eq:inprod}
(\bfa,\bfb)_\rho=\int_\calG \rho(\bfxi)\bfa(\bfxi)\cdot\bfb(\bfxi)\dxi
\end{equation}
with associated norm $\|\bfa\|=(\bfa,\bfa)_\rho^{1/2}$, where $\rho(\bfxi)\in\mathbb{R}$ is an arbitrary weight function. Consider all possible realization of $\bfu$, the associated mean vector of $\bfdelta$ is
\begin{equation}
\left\langle\bfdelta\right\rangle=\int_\calU \bfdelta(\bfxi,\bfu)p(\bfu)\du
\end{equation}
and the associated geometrical variance equals to
\begin{equation}\label{eq:sigma2}
\sigma^2=\left\langle\|\hat{\bfdelta}\|^2\right\rangle=\iint\limits_{\calU,\calG} \rho(\bfxi)\hat{\bfdelta}(\bfxi,\bfu)\cdot\hat{\bfdelta}(\bfxi,\bfu)p(\bfu)\dxi\du
\end{equation}
where $\hat{\bfdelta}=\bfdelta-\left\langle\bfdelta\right\rangle$, with $\left\langle\cdot\right\rangle$ the ensemble average \change{obtained integrating $\bfu$ over $\calU$}. 

The aim of KLE is to find an optimal basis of orthonormal functions for the linear representation of $\hat{\bfdelta}$:
\begin{equation}\label{eq:delta}
\hat{\bfdelta}(\bfxi,\bfu)\approx\sum_{k=1}^N x_k(\bfu)\bfphi_k(\bfxi)
\end{equation}
where 
\begin{equation}\label{eq:xred}
x_k(\bfu)=\left(\hat{\bfdelta},\bfphi_k\right)_\rho=\int_\calG \rho(\bfxi)\hat{\bfdelta}(\bfxi,\bfu)\cdot\bfphi_k(\bfxi)\dxi
\end{equation}
are the basis-function components usable as new (reduced) design variables. The optimality condition associated to the KLE refers to the geometric variance retained by the basis functions through Eq. \ref{eq:delta}. Combining Eqs. \ref{eq:sigma2}-\ref{eq:xred} yields
\begin{equation}
\sigma^2=\sum_{k=1}^\infty\left\langle\left(\hat{\bfdelta},\bfphi_k\right)_\rho^2\right\rangle
\end{equation}

The basis retaining the maximum variance is formed by those $\bfphi$, solutions of the variational problem
\begin{eqnarray}\label{eq:varp}
\underset{\bfphi\in\calL_\rho^2(\calG)}{\rm maximize} & \calJ(\bfphi)=\left\langle\left(\hat{\bfdelta},\bfphi_k\right)_\rho^2\right\rangle\\
\mathrm{subject \, to} & \left(\bfphi,\bfphi\right)_\rho^2=1\nonumber
\end{eqnarray}
which yields \cite{diez2015-CMAME}
\begin{equation}\label{eq:bie}
\calL\bfphi(\bfxi)=\int_\calG \rho(\bfxi')\left\langle\hat{\bfdelta}(\bfxi,\bfu)\otimes\hat{\bfdelta}(\bfxi',\bfu)\right\rangle\bfphi(\bfxi')\dxi'=\lambda\bfphi(\bfxi)
\end{equation}
where \change{$\otimes$ indicates the outer product and} $\calL$ is the self-adjoint integral operator whose eigensolutions define the optimal basis functions for the linear representation of Eq. \ref{eq:delta}. Therefore, its eigenfunctions (KL-modes) $\{\bfphi_k\}_{k=1}^\infty$ are orthonormal and form a complete basis for $\calL_\rho^2(\calG)$. Additionally, it may be proven that 
\begin{equation}
\sigma^2=\sum_{k=1}^\infty\lambda_k \qquad \mathrm{with} \qquad \lambda_k=\left\langle x_k^2\right\rangle
\end{equation}
where the eigenvalues $\lambda_k$ represent the variance retained by the associated basis function $\bfphi_k$, through its component $x_k$. Finally, the solution $\{\bfphi_k\}_{k=1}^\infty$ of Eq. \ref{eq:varp} are used to build a reduced dimensionality representation of the original design space; defining the desired confidence level $l$, with $0<l\leq1$, the number of reduced design variables $N$ in Eq. \ref{eq:delta} is selected such as
\begin{equation}\label{eq:lambda}
\sum_{k=1}^N \lambda_k\geq l\sum_{k=1}^\infty \lambda_k=l\sigma^2 \qquad \mathrm{with} \qquad \lambda_k\geq\lambda_{k+1}
\end{equation}

\subsection{Principal component analysis: dimensionality reduction in the discretized geometric space}\label{sec:pca}
Discretizing $\calG$ by $L$ elements \change{of measure $\Delta\calG_i$ (with $i=1,\dots,L$)}, 
sampling $\calU$ by a statistically convergent number of Monte Carlo (MC) realizations $S$, so that $\{\bfu_k\}_{k=1}^S  \sim p(\bfu)$, and organizing the discretization $\mathbf{d}(\bfxi,\bfu_k)$ of $\hat{\bfdelta}(\bfxi,\bfu_k)$ in a data matrix $\mathbf{D}$ of dimensionality $\left[3L\times S\right]$
\begin{equation}\label{eq:data}
\mathbf{D}=\left[
\begin{array}{ccc}
d_{1,\xi_1}(\bfu_1) & & d_{1,\xi_1}(\bfu_S)\\
\vdots & & \vdots\\
d_{L,\xi_1}(\bfu_1) & & d_{L,\xi_1}(\bfu_S)\\
d_{1,\xi_2}(\bfu_1) & & d_{1,\xi_2}(\bfu_S)\\
\vdots & \dots & \vdots\\
d_{L,\xi_2}(\bfu_1) & & d_{L,\xi_2}(\bfu_S)\\
d_{1,\xi_3}(\bfu_1) & & d_{1,\xi_3}(\bfu_S)\\
\vdots & & \vdots\\
d_{L,\xi_3}(\bfu_1) & & d_{L,\xi_3}(\bfu_S)\\
\end{array}
\right]
\end{equation}   
\change{where $d_{i,\xi_j}$ is the $j$-th component of the shape modification vector associated to the element $i$,} the integral problem of Eq. \ref{eq:varp} reduces to the generalized PCA of the data matrix $\mathbf{D}$, as follows
\begin{equation}\label{eq:pca}
\mathbf{A\change{G}WZ}=\mathbf{Z}\boldsymbol{\Lambda} \qquad \mathrm{with} \qquad \mathbf{A}=\frac{1}{S}\mathbf{DD}^\mathsf{T}
\end{equation} 
where $\mathbf{Z}$ and $\boldsymbol{\Lambda}$ are the eigenvectors and eigenvalues matrices of $\mathbf{A\change{G}W}$. \change{The discretization process of the problem \ref{eq:bie} and how it yields problem \ref{eq:pca} are discussed in \cite{diez2015-CMAME} and therefore not repeated here.}
\change{
The matrix $\mathbf{G}=\mathrm{diag}\left(\mathbf{G}_1,\mathbf{G}_2,\mathbf{G}_3 \right)$ is  block diagonal and has dimensionality $\left[3L\times 3L\right]$, with each $\left[L\times L\right]$ block ($k=1,2,3$) being a diagonal matrix itself
\begin{equation}
\mathbf{G}_k=\mathrm{diag}\left(\Delta\calG_1, \dots, \Delta\calG_L\right) 
\end{equation}
containing the measure $\Delta\calG_i$ of the $i$-th element.
}
Similarly, 
$\mathbf{W}=\mathrm{diag}\left(\mathbf{W}_1,\mathbf{W}_2,\mathbf{W}_3 \right)$ is a block diagonal matrix of dimensionality $\left[3L\times 3L\right]$, where each $\left[L\times L\right]$ block $\mathbf{W}_k$ ( $k=1,2,3$) is itself a diagonal matrix defined as
\begin{equation}\label{eq:weights}
\mathbf{W}_k=\mathrm{diag}\left({\rho_1}, \dots, {\rho_L}\right)
\end{equation}
with $\rho_i$  (for $i = 1,\dots, L$) the arbitrary weight given to each element.
 
\change{
It may be noted that if all elements of $\calG$ have equal measure $\Delta\calG$ and are assigned the same weight $\rho$, then the matrix $\mathbf{G}\mathbf{W}$ reduces to $\Delta\calG \, \rho \, \mathbf{I}$, where $\mathbf I$ is the identity matrix. Consequently, the problem \ref{eq:pca} reduces to the standard PCA of $\mathbf{D}$
\begin{equation}
\label{eq:pca_standard}
\mathbf{AZ}=\mathbf{Z}\boldsymbol{\Omega} \qquad \mathrm{with} \qquad \mathbf{\Omega}=\frac{1}{\Delta\calG \, \rho}\mathbf{\Lambda}
\end{equation}
}
\blue{Furthermore, the discretization of $\calG$ can be arbitrary and it does not necessarily have to coincide with that of the SDDO numerical solver. Nevertheless, the discretization used for the design-space dimensionality reduction can take in consideration the physics of the SDDO problems adding arbitrary weights to the elements, where most significant physical phenomena are involved.}

Finally, note that for both problems \ref{eq:pca} and \ref{eq:pca_standard} the matrix $\mathbf{Z}$ contains the discrete representation $\bfz_k$ of the desired eigenfunctions $\bfphi_k$.


\section{Parametric model embedding}\label{sec:III}

%
\begin{figure}[!t]
\centering
\includegraphics[width=0.9\textwidth]{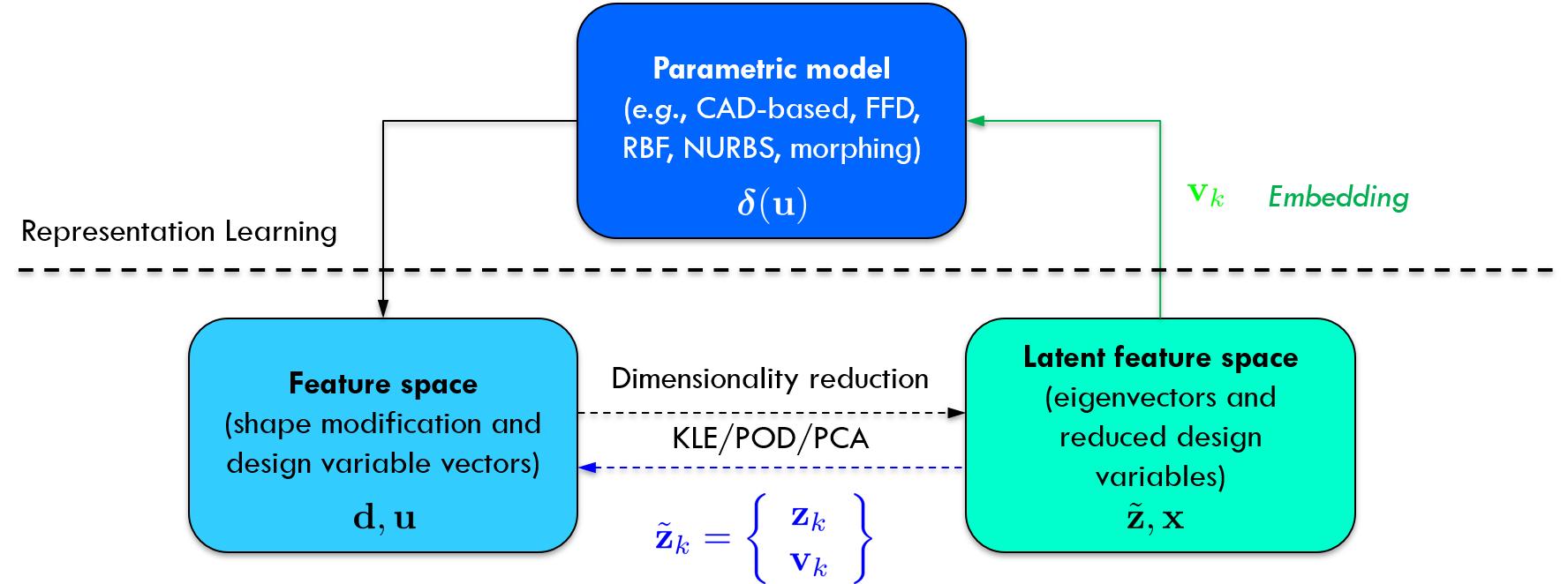}
\caption{Parametric model embedding concept}\label{fig:pme}
\end{figure}

In order to overcome KLE/PCA limitations, 
\change{starting from the discrete formulation presented in Section \ref{sec:pca},}
the data matrix $\mathbf{D}$ used for the dimensionality reduction procedure is augmented with the values of design parameters from the original parameterization \blue{(which is intrinsically finite)}, as conceptually shown in Fig. \ref{fig:pme}.

Practically, \change{defining $\hat{\bfu}=\bfu-\langle\bfu\rangle$ (as already done for the shape modification vector)}, the embedding is achieved \change{introducing} a new matrix $\mathbf{P}$ of dimensionality $\left[(3L+M)\times S\right]$ \change{as follows
\begin{equation}
\mathbf{P}=\left[
\begin{array}{c}
\mathbf{D}\\
\mathbf{U}
\end{array}
\right]
\qquad
\mathrm{with}
\qquad
\mathbf{U}=\left[
\begin{array}{ccc}
u_{1,1} & & u_{1,S}\\
\vdots & \cdots & \vdots\\
u_{M,1} & & u_{M,S}\\
\end{array}
\right]
\end{equation}   
}
where the matrix $\mathbf{U}$ of the original design variables is added to the data matrix $\mathbf{D}$ with a null weight \change{$\mathbf{W}_\bfu$ such that
\begin{equation}\label{eq:pme_weight}
\mathbf{W}_\bfu=\mathbf{0}
\qquad
\mathrm{and}
\qquad
\widetilde{\mathbf{W}}=\left[
\begin{array}{cc}
\mathbf{W} & \mathbf{0}\\
\mathbf{0} & \mathbf{W}_\bfu\\
\end{array}
\right]
\end{equation}   
}
and so recasting Eq. \ref{eq:pca} to
\change{
\begin{equation}\label{eq:pme_pca}
\widetilde{\mathbf{A}}\widetilde{\mathbf{G}}\widetilde{\mathbf{W}}\widetilde{\mathbf{Z}}=\widetilde{\mathbf{Z}}\widetilde{\boldsymbol{\Lambda}} \qquad \mathrm{with} \qquad \widetilde{\mathbf{A}}=\frac{1}{S}\mathbf{PP}^\mathsf{T} 
\end{equation} 
where
\begin{equation}\label{eq:pme}
\widetilde{\mathbf{G}}=\left[
\begin{array}{cc}
\mathbf{G} & \mathbf{0}\\
\mathbf{0} & \mathbf{I}\\
\end{array}
\right]
\qquad
\mathrm{and}
\qquad
\widetilde{\mathbf{Z}} = \left[\tilde{\bfz}_1 \,\,\, \dots \,\,\, \tilde{\bfz}_S \right] \qquad \mathrm{with} \qquad 
\tilde{\bfz}_k=\left\{
\begin{array}{c}
\bfz_k\\
\bfv_k
\end{array}
\right\}
\end{equation} 
}

\change{Having given a null weight to $\mathbf{U}$ does not remove the contribution of the design variables from the inner product, but just cancels as many columns as $M$ from the  matrix $\widetilde{\mathbf{A}}\widetilde{\mathbf{G}}\widetilde{\mathbf{W}}$, thus Eqs. \ref{eq:pca} and \ref{eq:pme} provides the same eigenvalues ($\boldsymbol{\Lambda}=\widetilde{\boldsymbol{\Lambda}}$) and geometrical components of the eigenvectors ($\bfz_k$, except for a multiplicative constant). The proof of the equivalence between KLE and PME is given in \ref{sec:appA}.} In addition and as desired, the solution of Eq. \ref{eq:pme} provides the eigenvector components $\bfv_k$ that embeds the original design variables $\bfu$. \blue{It may be emphasized that the choice of giving zero weight to the design variables stems from the fact that we are not interested in their variance per se (prior distribution), but to the resulting variance of the shape modification vector (posterior distribution), which is at the basis of the current dimensionality reduction method.}

In order to reconstruct at least all the samples in $\mathbf{D}$, the reduced design variables $\mathbf{x}$ are bounded such as $\mathrm{inf}\{\boldsymbol{\alpha}_k\}\leq \bfx \leq \mathrm{sup}\{\boldsymbol{\alpha}_k\}$, with
\begin{equation}
\boldsymbol{\alpha}_k=\mathbf{p}_k^\mathsf{T}\widetilde{\mathbf{G}}\widetilde{\mathbf{W}}\widetilde{\mathbf{Z}}' \qquad \mathrm{for} \,\,\, k=1,\dots,S
\end{equation}
where $\widetilde{\mathbf{Z}}'$ contains only the first $N$ eigenvectors of $\widetilde{\mathbf{Z}}$, retaining the desired level of variance of the original design space, and $\mathbf{p}_k$ is a column of the matrix $\mathbf{P}$. 

The PME of the original design variables 
is finally achieved by
\change{
\begin{equation}\label{eq:recu}
\bfu \approx \hat{\bfu}=\langle{\bfu}\rangle+\sum_{k=1}^N x_k\bfv_k
\end{equation}
\blue{where the eigenvectors component $\bfv_k$ embeds (or contains) the reduced-order representation of the original design parameterization.}
It may be noted that the overall methodology is independent from the specific shape modification method. This is seen as a black box by PME.
}
\section{Example application\blue{s}}\label{sec:IV}
The PME method is demonstrated, \blue{as a proof of concept,} for the shape reparameterization and optimization of the \blue{NACA 0012 airfoil and the DTMB 5415 model (an open-to-public naval combatant), both widely used as benchmarks in the aerodynamic (see, e.g., \cite{he2019robust}) and the ship hydrodynamic (see, e.g.,  \cite{grigoropoulos2017-MARINE}) communities}. 

\subsection{\blue{Shape optimization problems formulation}}
The PME effectiveness is shown for the solution of \blue{two} SDDO \blue{problems} formulated \blue{for the drag coefficient reduction in subsonic condition and the resistance reduction in calm water of NACA 0012 and DTMB 5415, respectively. Having no prior knowledge of design space and their performance, a uniform PDF for $\bfu$ is used for both SDDO problems, whose details are provided in the following.} 

\subsubsection{\blue{NACA 0012 problem}}
\begin{figure}[!h]
\centering
\includegraphics[width=0.75\textwidth]{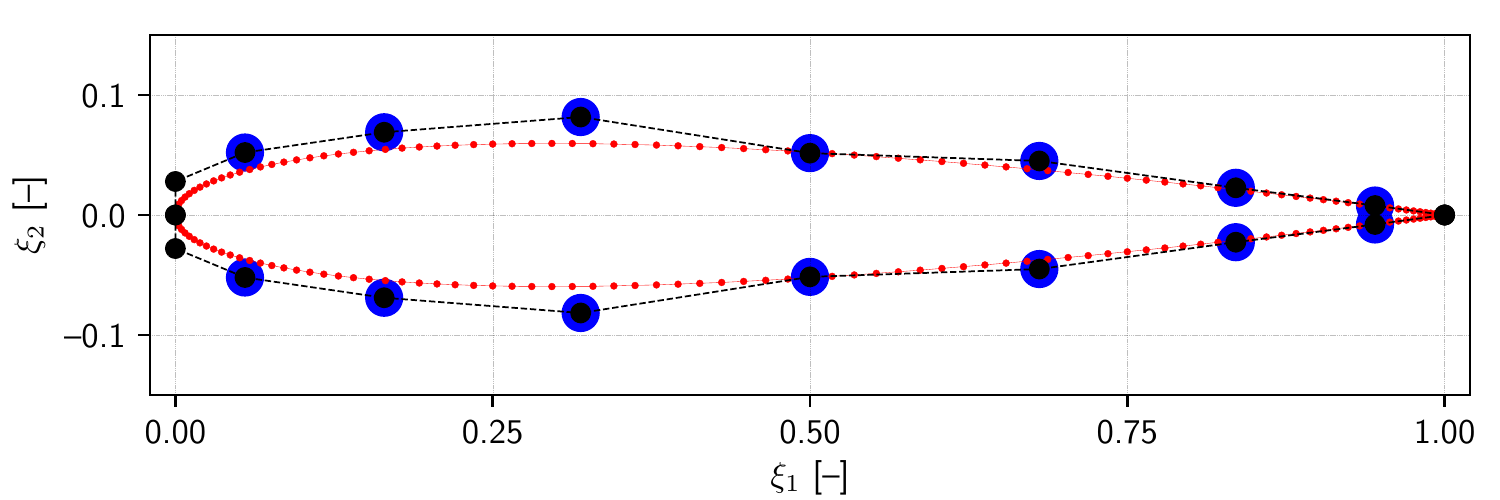}
\caption{\blue{Bezier-curves design space definition and grid nodes used for the dimensionality reduction}}\label{fig:bezier}
\end{figure}

\blue{The objective is the reduction of the drag coefficient $C_D$ at fixed Mach number and lift coefficient ($C_L$), respectively equal to 0.3 and 0.5. The optimization problem reads}
\blue{
\begin{equation}\label{eq:SDDOprob_foil}
\begin{array}{rll}
\mathrm{minimize}      & C_D(\bfu)\\ 
\mathrm{subject \, to} & C_L(\bfu) = 0.5\\
\mathrm{and \, to}     & -0.1\leq C_M(\bfu)\leq 0\\ 
                       & \bfu^l\leq\bfu\leq\bfu^u,
\end{array}
\end{equation}  
}
\blue{where $C_M$ is the pitching moment coefficient and $\bfu^l$ and $\bfu^u$ are the lower and upper bounds of $\bfu$, respectively. The original design space is defined by $M=14$ design variables, based on Bezier curves of degree ten \cite{mohebbi2014aerodynamic}. Specifically, the airfoil coordinates are defined by two Bezier curves, one for the suction and one for pressure sides. Each of these curves has 11 control points, shown in black in Fig. \ref{fig:bezier}, but only 7 are active (see blue point in Fig. \ref{fig:bezier}) with only $\xi_2$ degree of freedom (DoF); their lower and upper bounds are -0.9 and 0.9.}

\blue{The objective function is evaluated by XFOIL \cite{Drela1989},
a boundary-layer (BL) solver with viscous correction. Specifically, the boundary layers and wake are described with a two-equation lagged dissipation integral BL formulation and an envelope $e^n$ transition criterion. Simulations are performed with Reynolds number equal to $5E6$ at fixed $C_L$, so the solver iteratively adjust the angle of attack $\alpha$ to match the desired $C_L$. The airfoil is discretized by 182 nodes (see red circle in Fig. \ref{fig:bezier}.}
\begin{figure}[!t]
\centering
\includegraphics[width=0.75\textwidth]{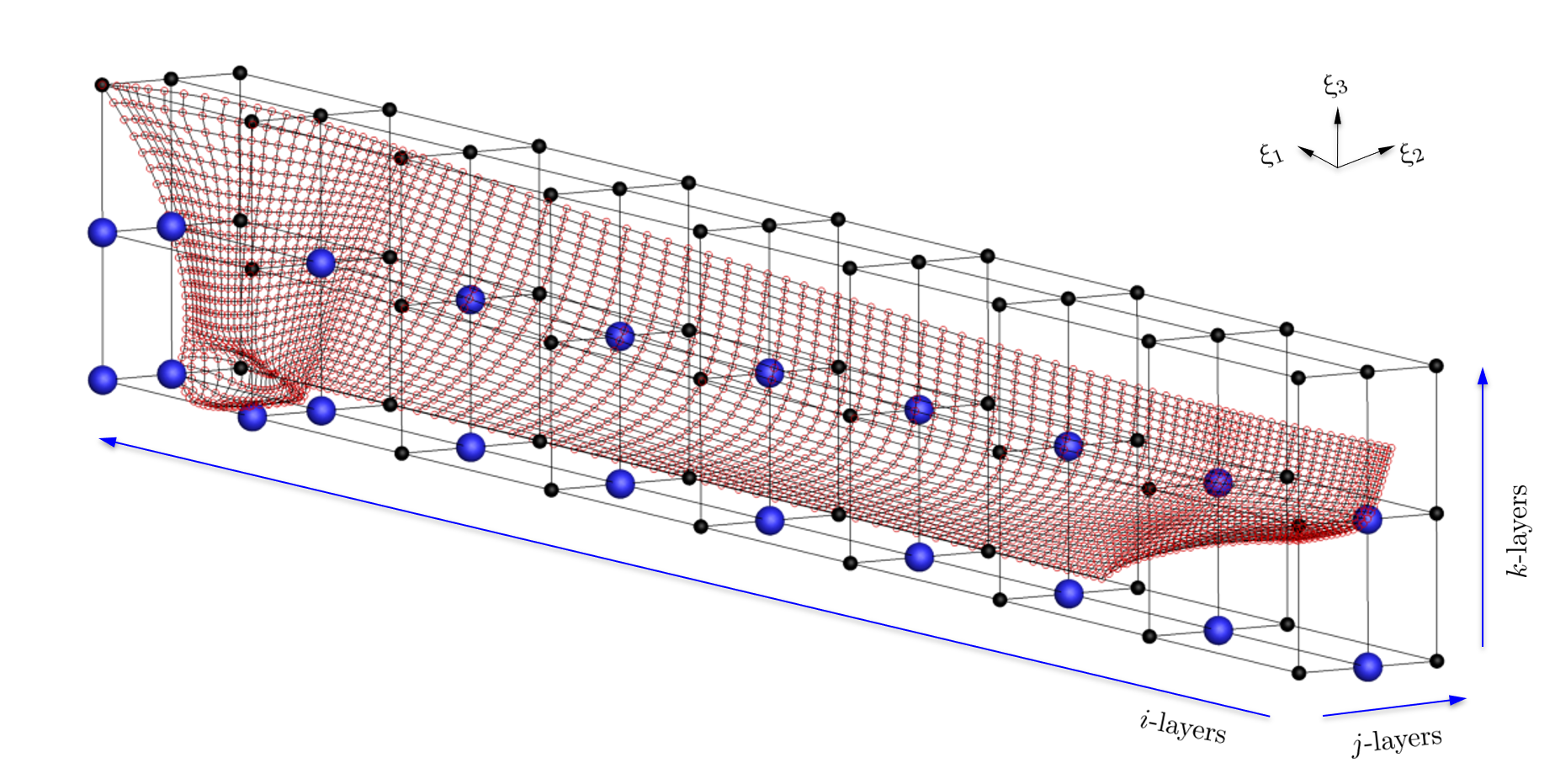}
\caption{FFD design spaces definition and grid nodes used for the dimensionality reduction}\label{fig:ffd}
\end{figure}
\begin{table}[!t]
\centering
\caption{FFD design variables definition}\label{tab:ffd}
\begin{tabular}{ccccccc}
\toprule
Design variable & $i$-layer & $j$-layer & $k$-layer & DoF & Lower bound & Upper bound\\
\midrule
$u_{1}$	&1	&2	&1	&$\xi_2$	&-0.500	&0.500\\
$u_{2}$	&2	&2	&1	&$\xi_2$	&-0.500	&0.500\\
$u_{3}$	&3	&2	&1	&$\xi_2$	&-0.500	&0.500\\
$u_{4}$	&4	&2	&1	&$\xi_2$	&-0.500	&0.500\\
$u_{5}$	&5	&2	&1	&$\xi_2$	&-0.500	&0.500\\
$u_{6}$	&6	&2	&1	&$\xi_2$	&-0.500	&0.500\\
$u_{7}$	&7	&2	&1	&$\xi_2$	&-0.500	&0.500\\
$u_{8}$	&8	&2	&1	&$\xi_2$	&-0.500	&0.500\\
$u_{9}$	&9	&2	&1	&$\xi_2$	&-0.500	&0.500\\
$u_{10}$&1	&2	&2	&$\xi_2$	&-0.500	&0.500\\
$u_{11}$&2	&2	&2	&$\xi_2$	&-0.500	&0.500\\
$u_{12}$&3	&2	&2	&$\xi_2$	&-0.500	&0.500\\
$u_{13}$&4	&2	&2	&$\xi_2$	&-0.500	&0.500\\
$u_{14}$&5	&2	&2	&$\xi_2$	&-0.500	&0.500\\
$u_{15}$&6	&2	&2	&$\xi_2$	&-0.500	&0.500\\
$u_{16}$&7	&2	&2	&$\xi_2$	&-0.500	&0.500\\
$u_{17}$&8	&2	&2	&$\xi_2$	&-0.500	&0.500\\
$u_{18}$&9	&2	&2	&$\xi_2$	&-0.500	&0.500\\
$u_{19}$&9	&1	&2	&$\xi_3$	&-0.250	&0.250\\
$u_{20}$&9	&1	&1	&$\xi_1$	&-0.025	&0.025\\
$u_{21}$&9	&1	&1	&$\xi_3$	&-0.100	&0.100\\
$u_{22}$&8	&1	&1	&$\xi_1$	&-0.025	&0.025\\
\bottomrule
\end{tabular}
\end{table}

\subsubsection{\blue{DTMB 5415 problem}}
\blue{The optimization pertains to the resistance $R$ reduction in calm water at Froude equal to 0.28 and reads}
\begin{equation}\label{eq:SDDOprob}
\begin{array}{rll}
\mathrm{minimize}      & R(\bfu)\\ 
\mathrm{subject \, to} & L_{\rm pp} = L_{\rm pp,0}\\
\mathrm{and \, to}     & \nabla(\bfu)\geq\nabla_0\\ 
                       & |\Delta B(\bfu)| \leq 5\%B_0\\ 
                       & |\Delta T(\bfu)| \leq 5\%T_0\\ 
                       & V(\bfu) \geq V_0, \\
                       & \bfu^l\leq\bfu\leq\bfu^u,
\end{array}
\end{equation}  
where $L_{\rm pp}$ is the length between \change{perpendiculars}, $\nabla$ the \change{volume displaced}, $B$ the overall beam, $T$ the \change{draught}, and $V$ the volume reserved for the sonar in the bow dome; finally, $\bfu^l$ and $\bfu^u$ are the lower and upper bounds of $\bfu$, respectively.
 
The design space has been defined within the activities of the NATO Science and Technology Organization, Applied Vehicle Technology (AVT), Research Task Group (RTG) 331 on “Goal-Driven, Multi-Fidelity Approaches for Military Vehicle System-Level Design'' \cite{beran2020-AIAA}. The original design space is formed by $M = 22$ design variables, defined by the FFD method \cite{sederberg1986-CGIT}. Specifically, the demi-hull is put in a lattice of $9\times 3 \times 3$ nodes in the $\xi_1\xi_2\xi_3$ reference system, Fig. \ref{fig:ffd}. Note that the FFD lattice perfectly fit the demi-hull maximum dimension. Only 21 nodes are active (see blue sphere in Fig. \ref{fig:ffd}) and their degrees of freedom (DoF), associated to the design variables, are summarized in Tab. \ref{tab:ffd} (only one active node has two DoF, all the others have one DoF). Note that the lower and upper bounds for $\bfu$ are provided normalizing the lattice in a unit cube. \blue{It may be noted that the FFD lattice definition has been introduced to distribute and replicate the shape modifications of a non-free CAD parametric model \cite{grigoropoulos2017-MARINE}, already use within the NATO-AVT-204 on “Assess the Ability to Optimize Hull Forms of Sea Vehicles for Best Performance in a Sea Environment''.}

The objective function is evaluated by a potential flow solver, developed at CNR-INM, whose details can be found in \cite{bassanini1994-SMI}. Specifically, using Dawson (double-model) linearization \cite{dawson1977-NSH}, the resistance of the ship is evaluated as the sum of wave and friction components: the former is evaluated by a standard pressure integral, whereas the latter is estimated using a flat-plate approximation, based on the local Reynolds number. Simulations are performed at even keel condition for the demi-hull only, taking advantage of symmetry about the $\xi_1\xi_3$-plane. The computational domain for the free-surface is defined within 0.5$L_{\rm pp}$ upstream, 1.5$L_{\rm pp}$ downstream, and 1$L_{\rm pp}$ sideways; $90\times 25$ grid nodes (see red circle in Fig. \ref{fig:ffd}) are used for the hull grid, whereas $75\times 22$ nodes are used for the free-surface discretization.

\begin{figure}[!b]
\centering
\includegraphics[width=0.219\textwidth]{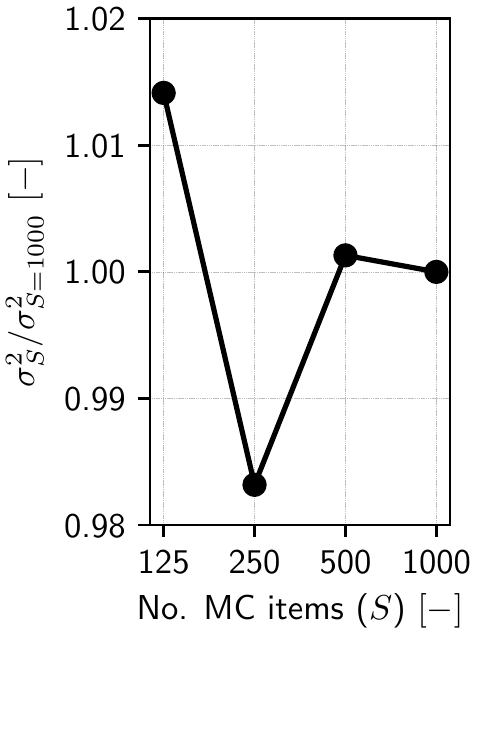}
\includegraphics[width=0.439\textwidth]{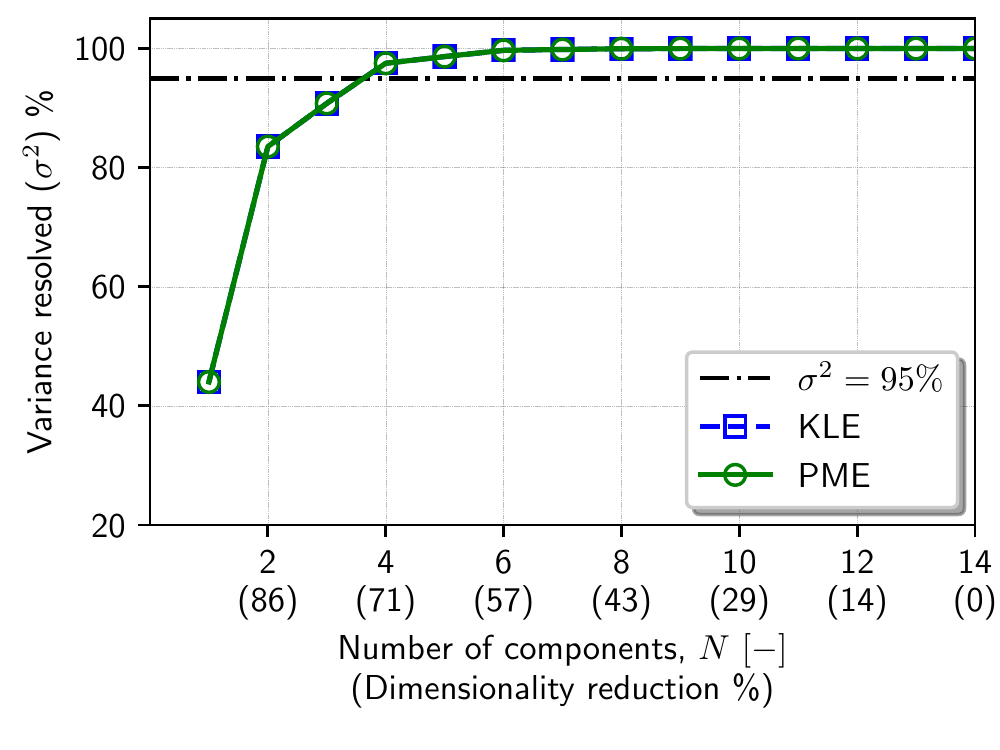}
\includegraphics[width=0.329\textwidth]{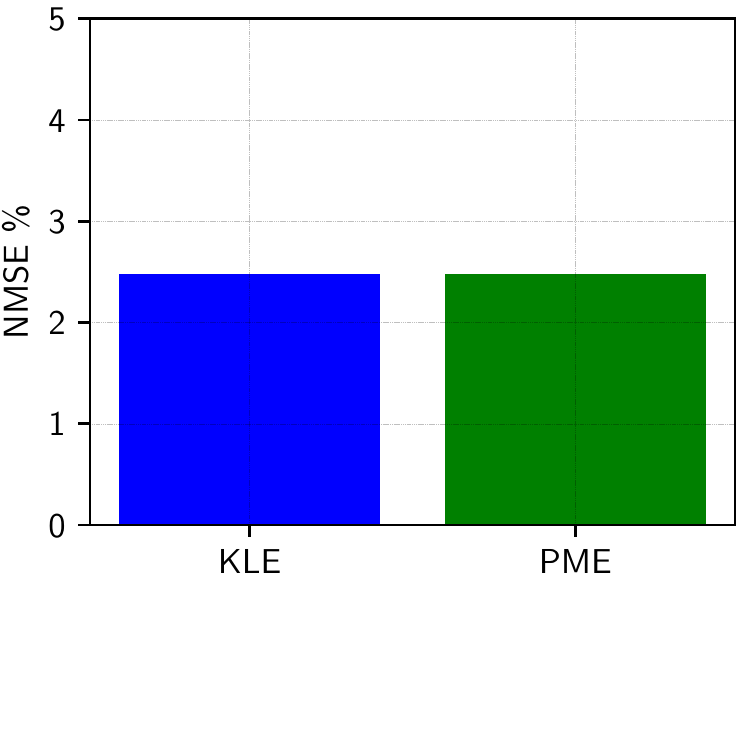}\\
\includegraphics[width=0.219\textwidth]{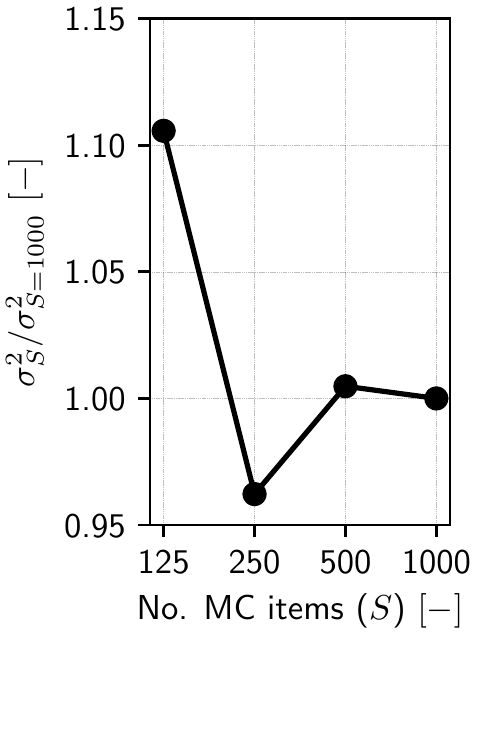}
\includegraphics[width=0.439\textwidth]{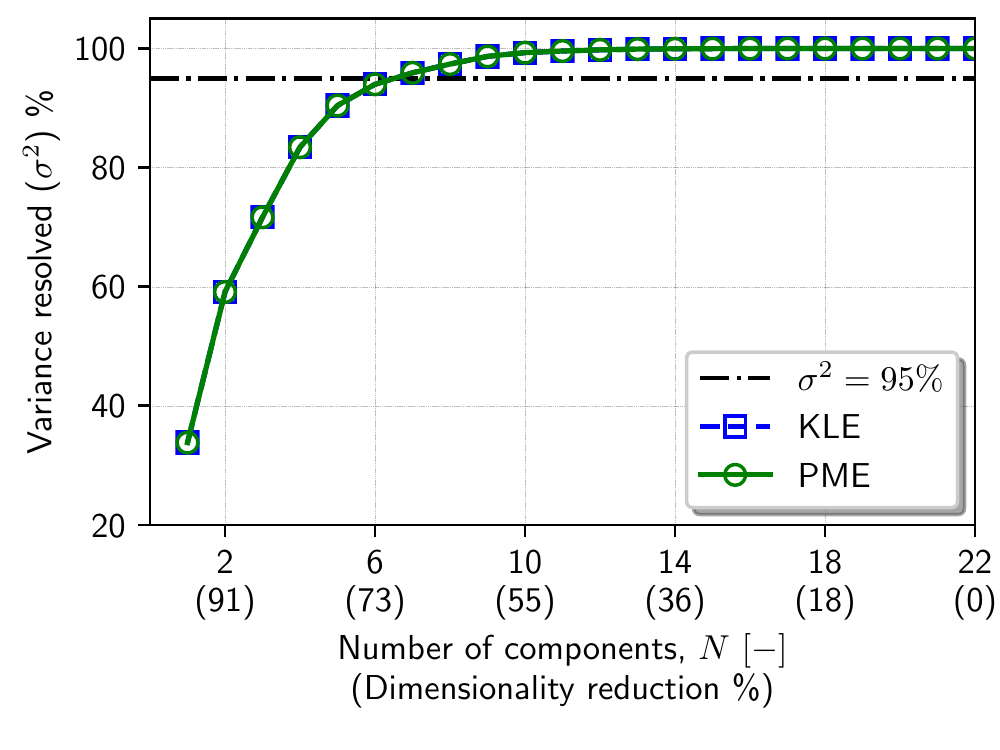}
\includegraphics[width=0.329\textwidth]{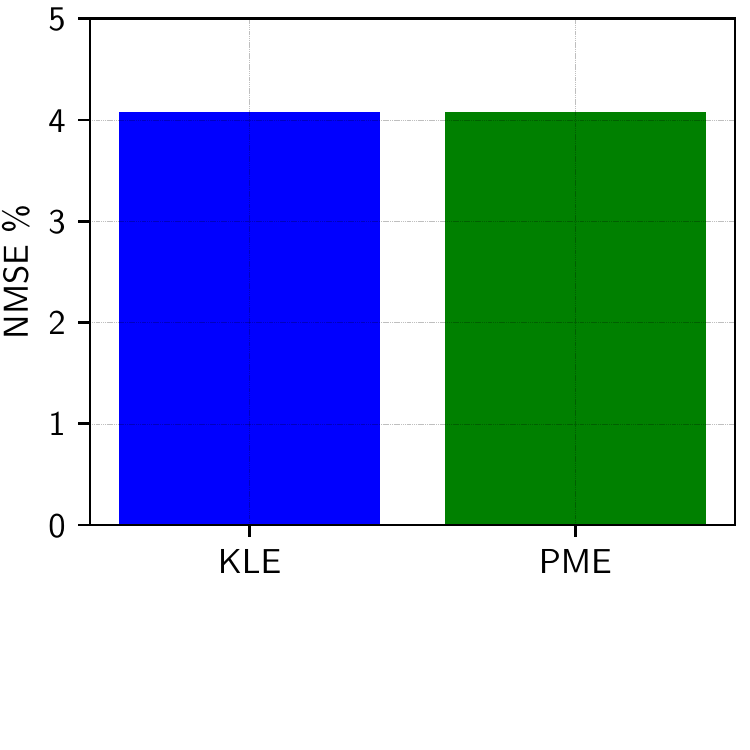}\\
\caption{\blue{Design-space dimensionality reduction results for NACA 0012 (top) and DTMB 5415 (bottom): }geometrical variance convergence towards the number of MC items (left), variance resolved as a function of the number of reduced design variables (center) and NMSE (right) obtained with KLE and PME}\label{fig:eigsum}
\end{figure}

\subsection{Dimensionality reduction results}
The PCA is trained by a set of $S = 1000$ MC items, following Diez and Serani \cite{diez2021-UQOP}, where a parametric and statistical analysis conditional to the number of MC samples has been conducted. \change{The convergence of the geometrical variance towards the number of MC items is provided in Fig. \ref{fig:eigsum} (left)}. For the sake of simplicity, in the following, KLE refers to the dimensionality reduction based on the shape modification vector only, whereas PME refers to the dimensionality reduction using both the shape modification and the design variables vectors. 
\blue{For the sake of simplicity, the data matrices $\mathbf{D}$ (see Eq. \ref{eq:data}), for both design spaces, are based on the geometry discretization used for the numerical solver. Specifically, $\mathbf{D}$ has dimensions $\left[364\times 1000\right]$ and $\left[6750\times 1000\right]$, for NACA 0012 and DTMB 5415, respectively. Furthermore, for the DTMB 5415, a} weight coefficient $\rho_i=1$ (see Eq. \ref{eq:weights}) is imposed for all the grid nodes below the water line, while a null weight ($\rho_i=0$) is used for the nodes above. 
\begin{figure}[!t]
\centering
\includegraphics[width=0.45\textwidth]{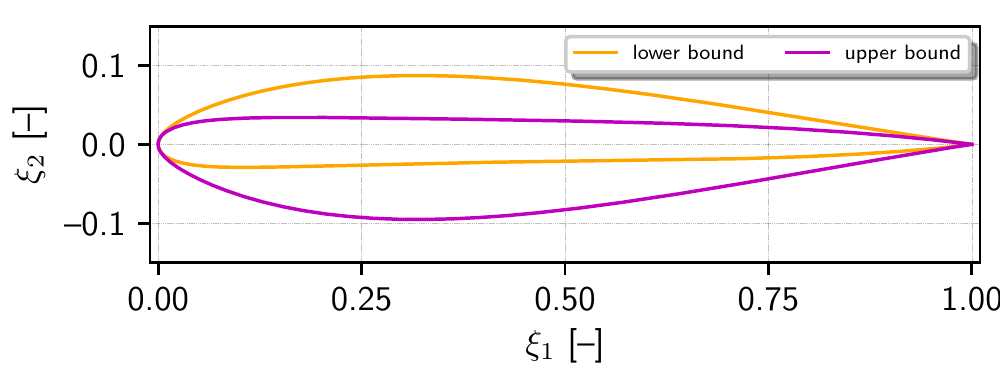}
\includegraphics[width=0.45\textwidth]{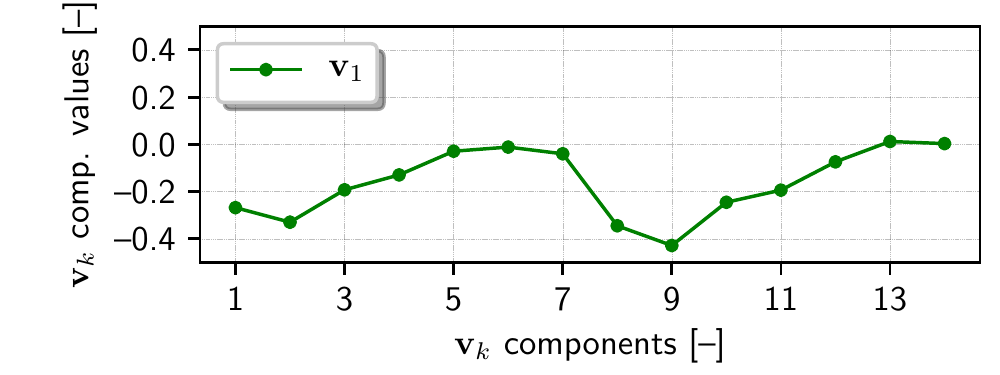}\\
\includegraphics[width=0.45\textwidth]{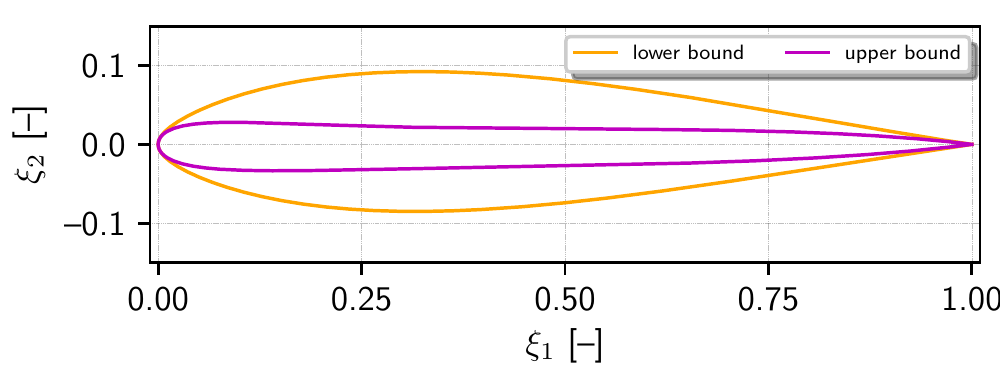}
\includegraphics[width=0.45\textwidth]{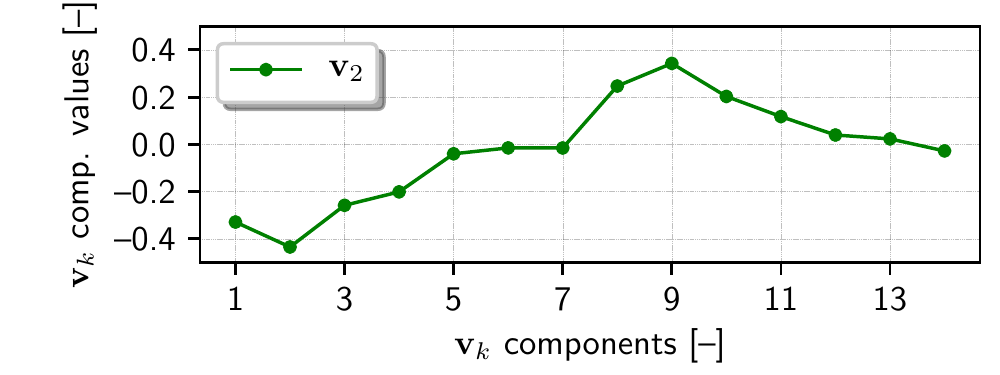}\\
\includegraphics[width=0.45\textwidth]{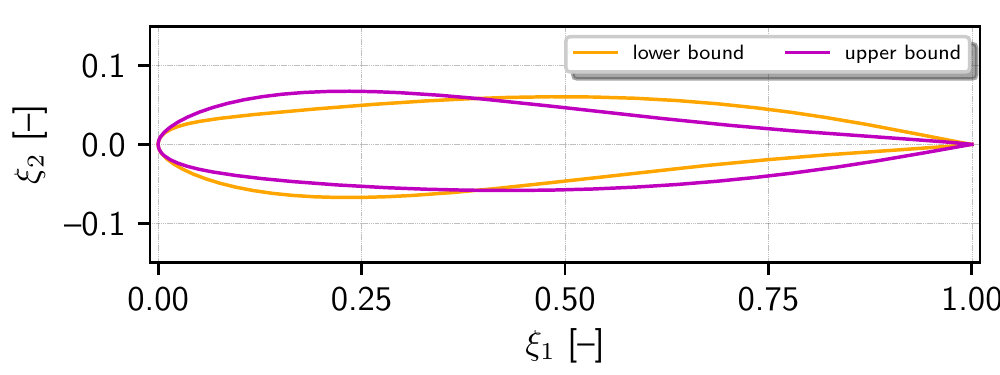}
\includegraphics[width=0.45\textwidth]{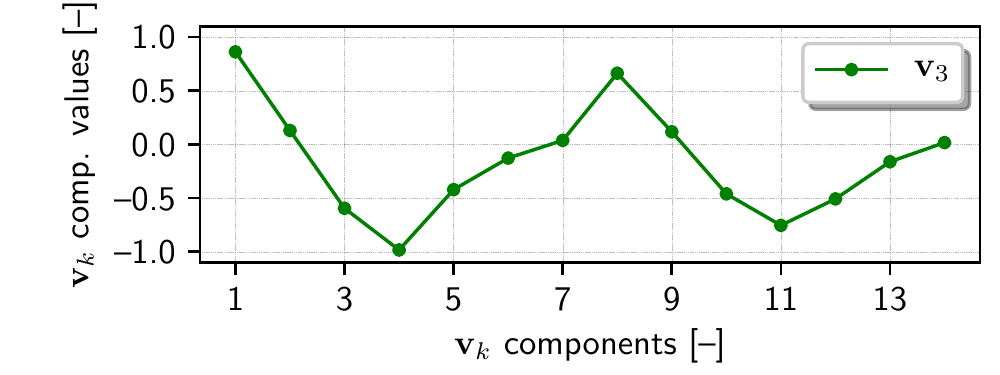}\\
\includegraphics[width=0.45\textwidth]{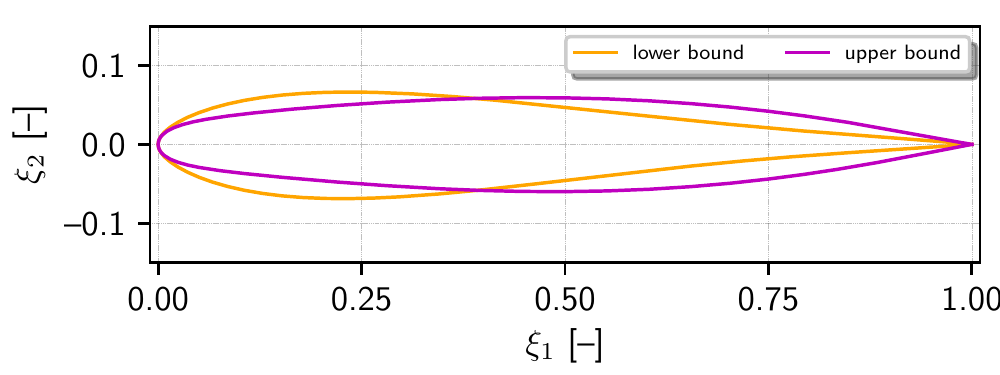}
\includegraphics[width=0.45\textwidth]{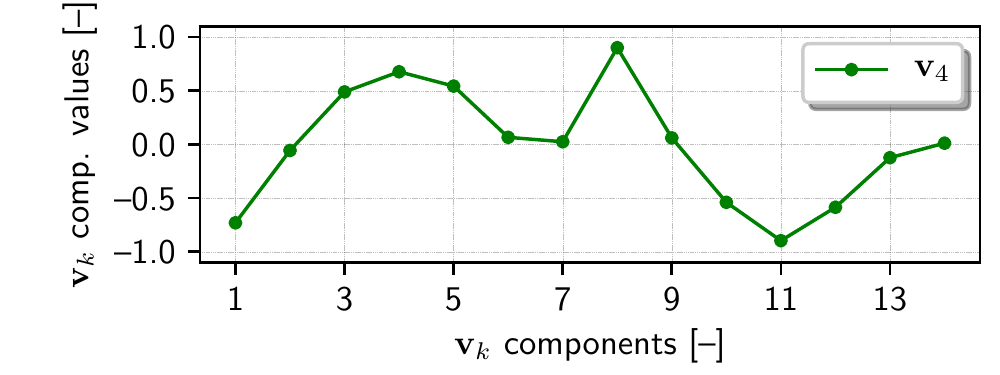}
\caption{\blue{Design-space dimensionality reduction modes for NACA 0012: (left) the shape modification vector modes $\mathbf{z}_k$ and (right) modes $\mathbf{v}_k$ that embeds the original design variables, for $k=1,\dots,4$ (from top to bottom)}}\label{fig:eigg_foil}
\end{figure}
\begin{figure}[!t]
\centering
\includegraphics[width=0.45\textwidth]{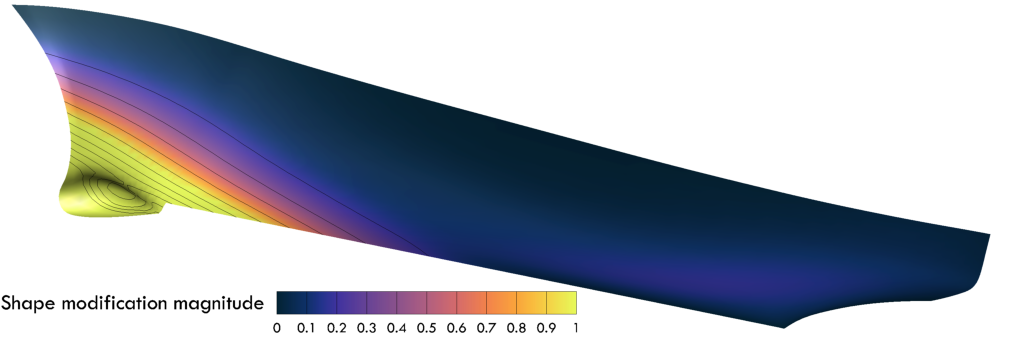}
\includegraphics[width=0.45\textwidth]{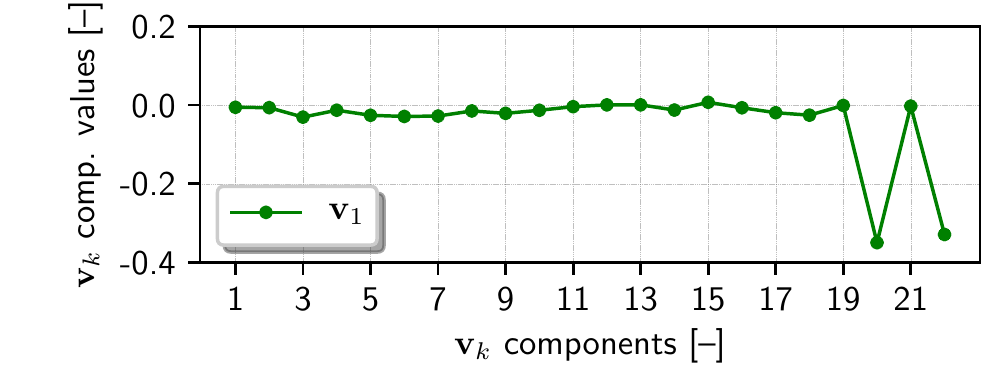}\\
\includegraphics[width=0.45\textwidth]{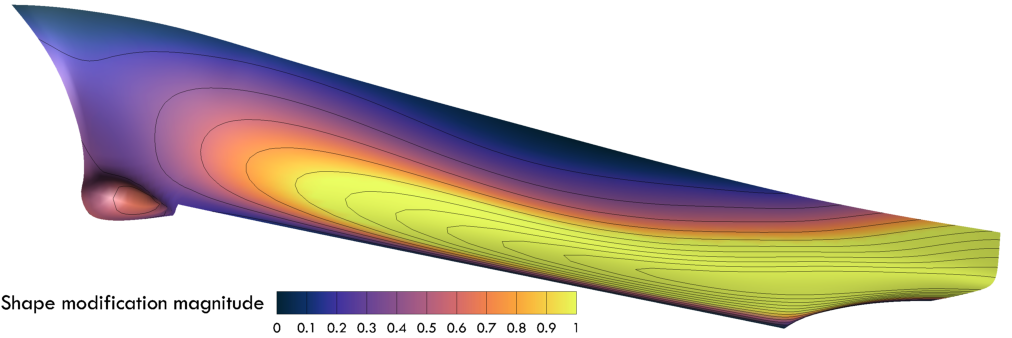}
\includegraphics[width=0.45\textwidth]{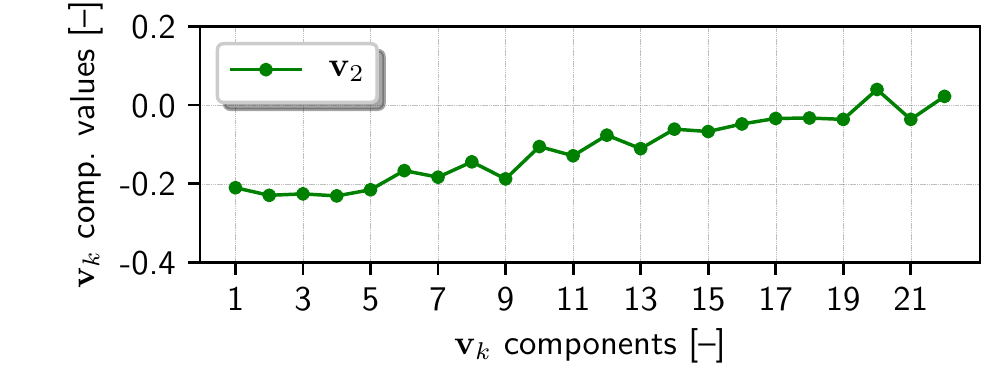}\\
\includegraphics[width=0.45\textwidth]{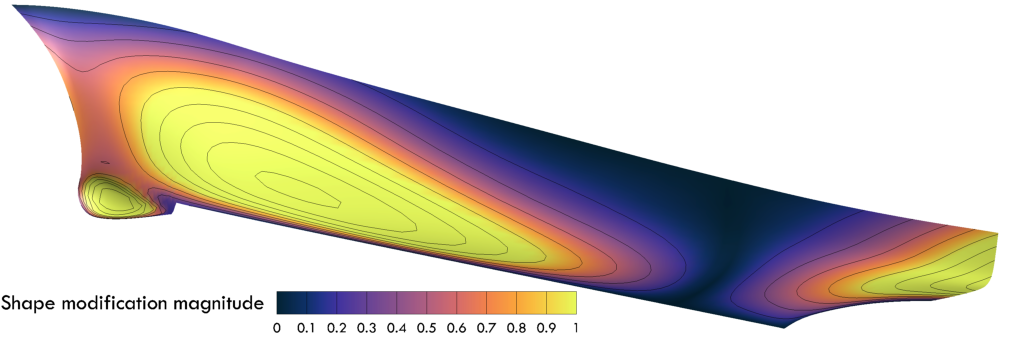}
\includegraphics[width=0.45\textwidth]{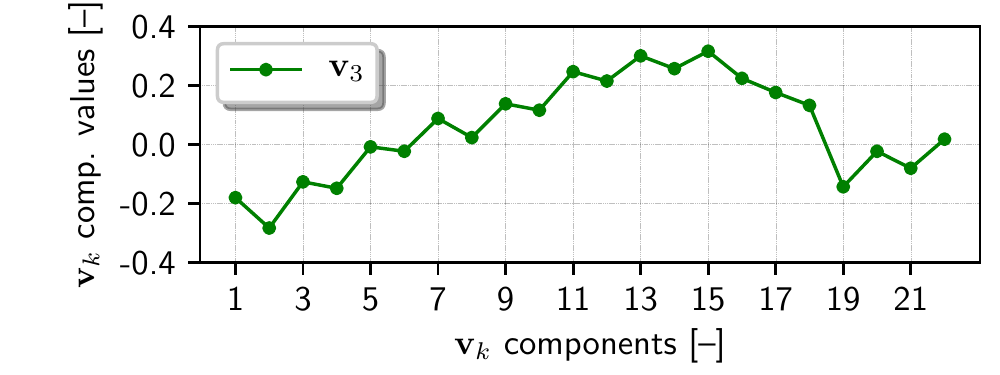}\\
\includegraphics[width=0.45\textwidth]{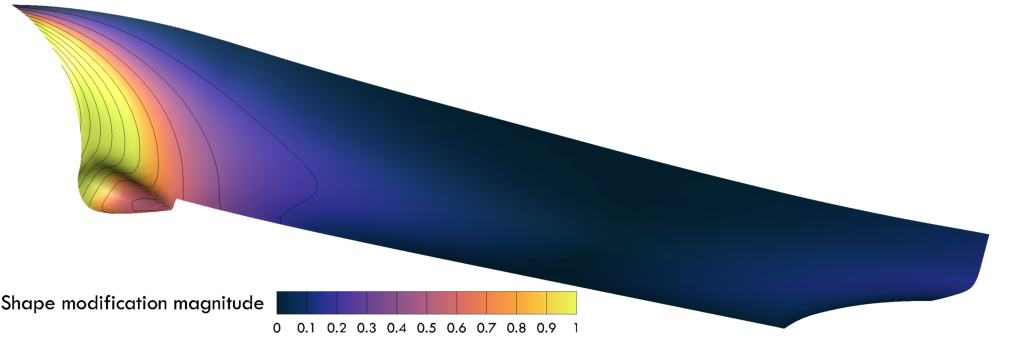}
\includegraphics[width=0.45\textwidth]{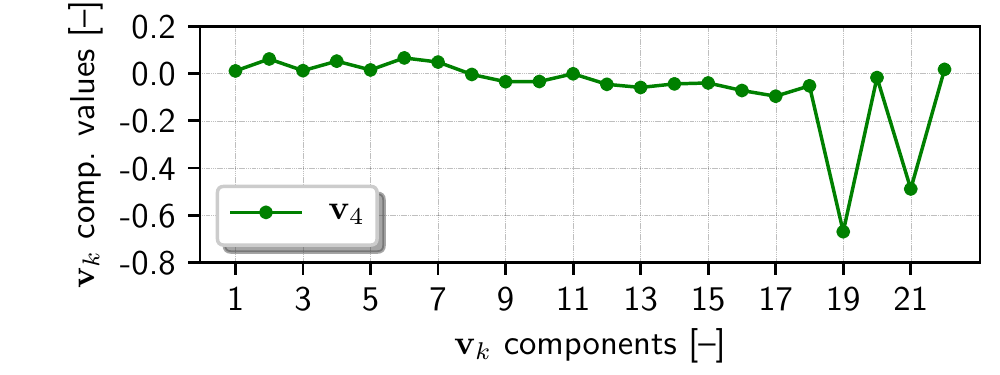}\\
\includegraphics[width=0.45\textwidth]{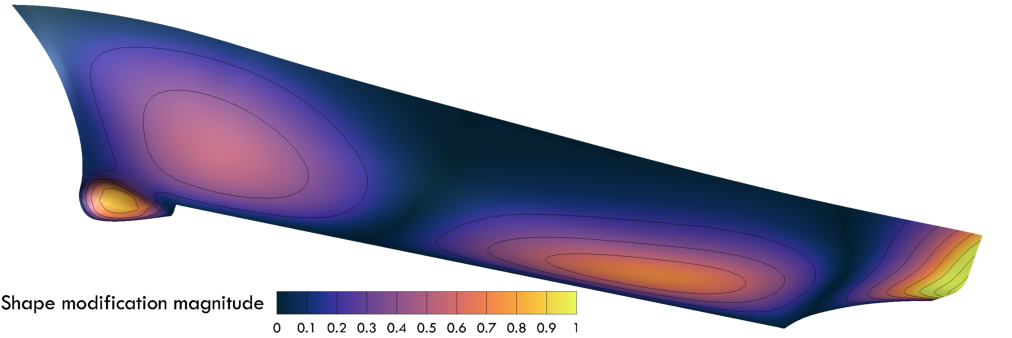}
\includegraphics[width=0.45\textwidth]{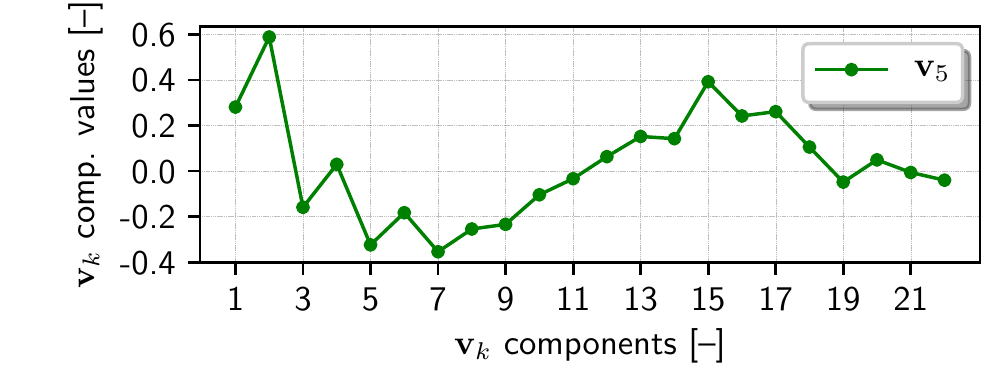}\\
\includegraphics[width=0.45\textwidth]{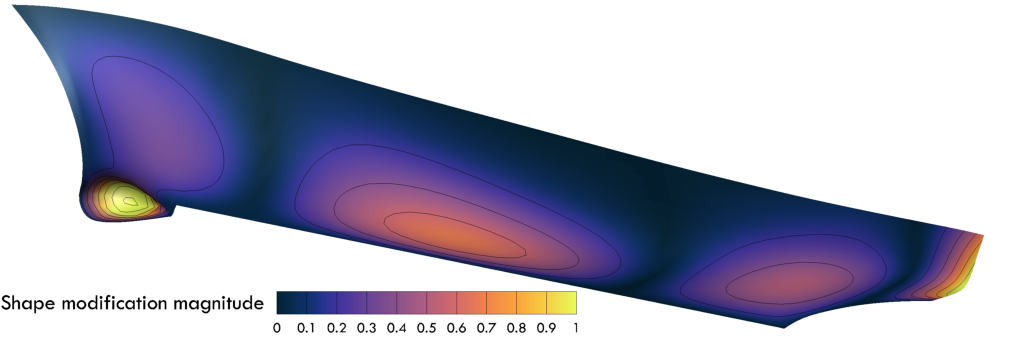}
\includegraphics[width=0.45\textwidth]{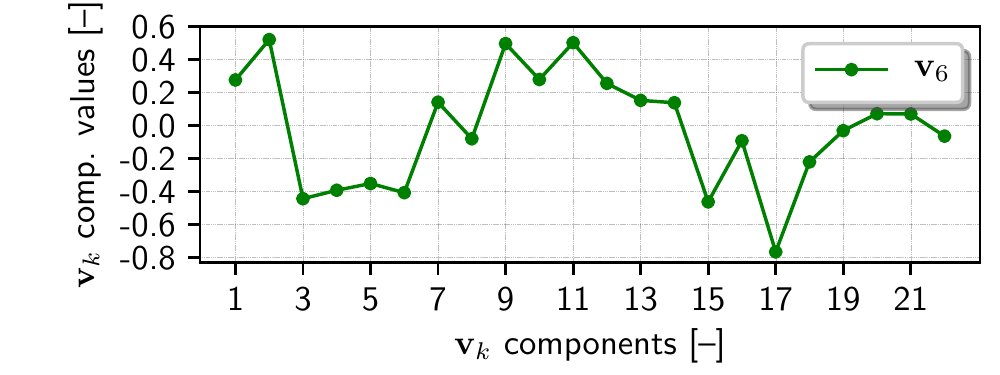}\\
\includegraphics[width=0.45\textwidth]{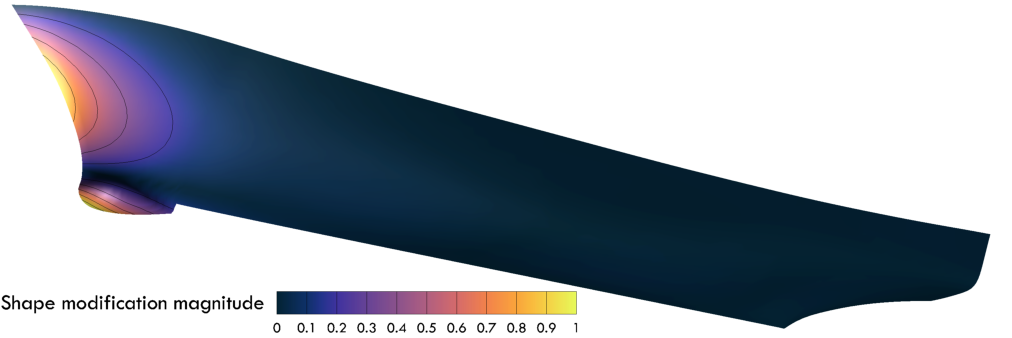}
\includegraphics[width=0.45\textwidth]{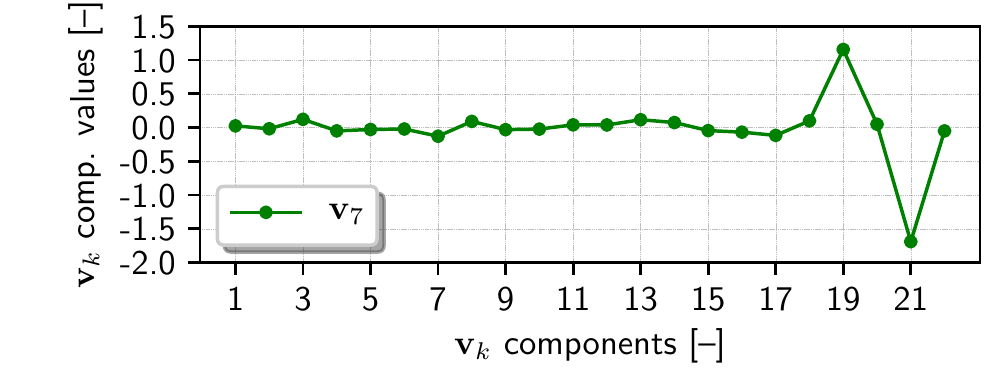}\\
\caption{Design-space dimensionality reduction modes \blue{for DTMB 5415}: (left) the shape modification vector modes $\mathbf{z}_k$ and (right) modes $\mathbf{v}_k$ that embeds the original design variables, for $k=1,\dots,7$ (from top to bottom)}\label{fig:eigg}
\end{figure}
\begin{figure}[!t]
\centering
\includegraphics[width=0.45\textwidth]{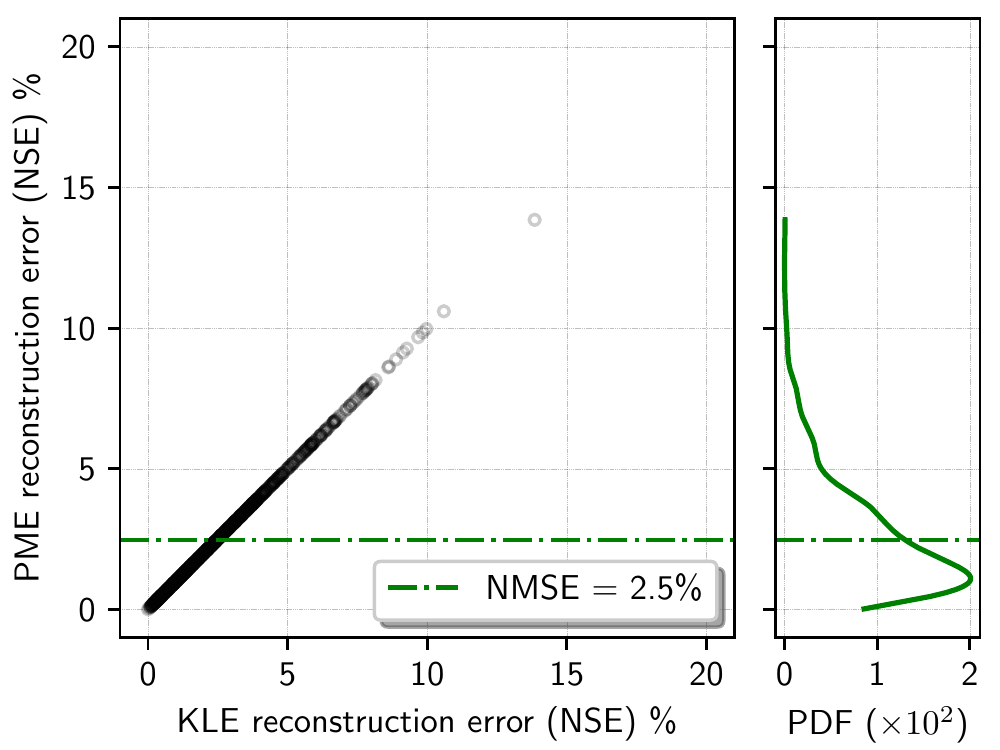}
\includegraphics[width=0.45\textwidth]{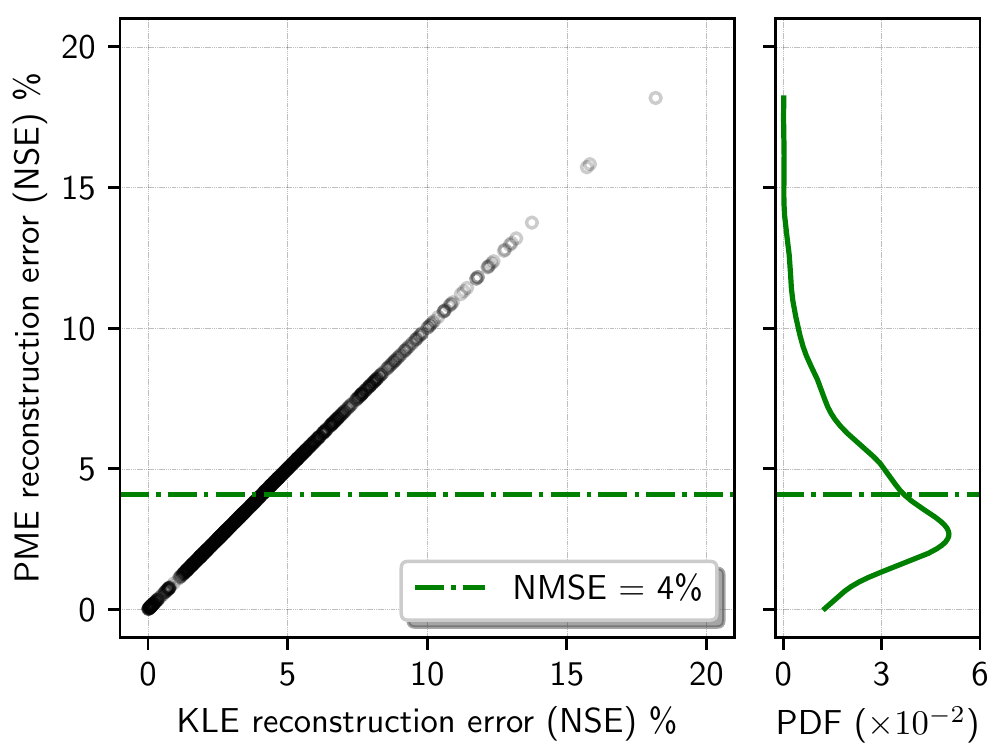}
\caption{Correlation of KLE and PME reconstruction errors for geometry data set $\mathbf{D}$ \blue{of NACA 0012 (left) and DTMB 5415 (right)}, along with PDF}\label{fig:corr}
\end{figure}
\begin{figure}[!t]
\centering
\includegraphics[width=0.32\textwidth]{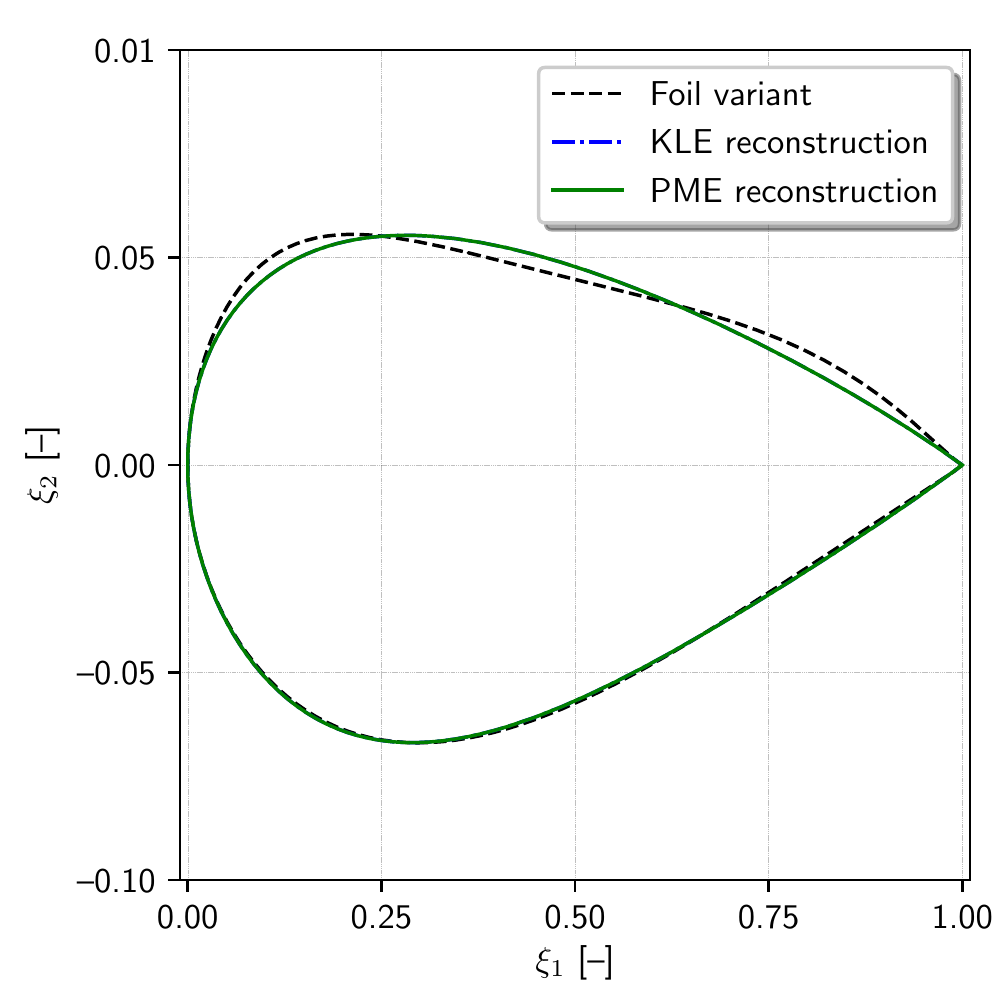}
\includegraphics[width=0.32\textwidth]{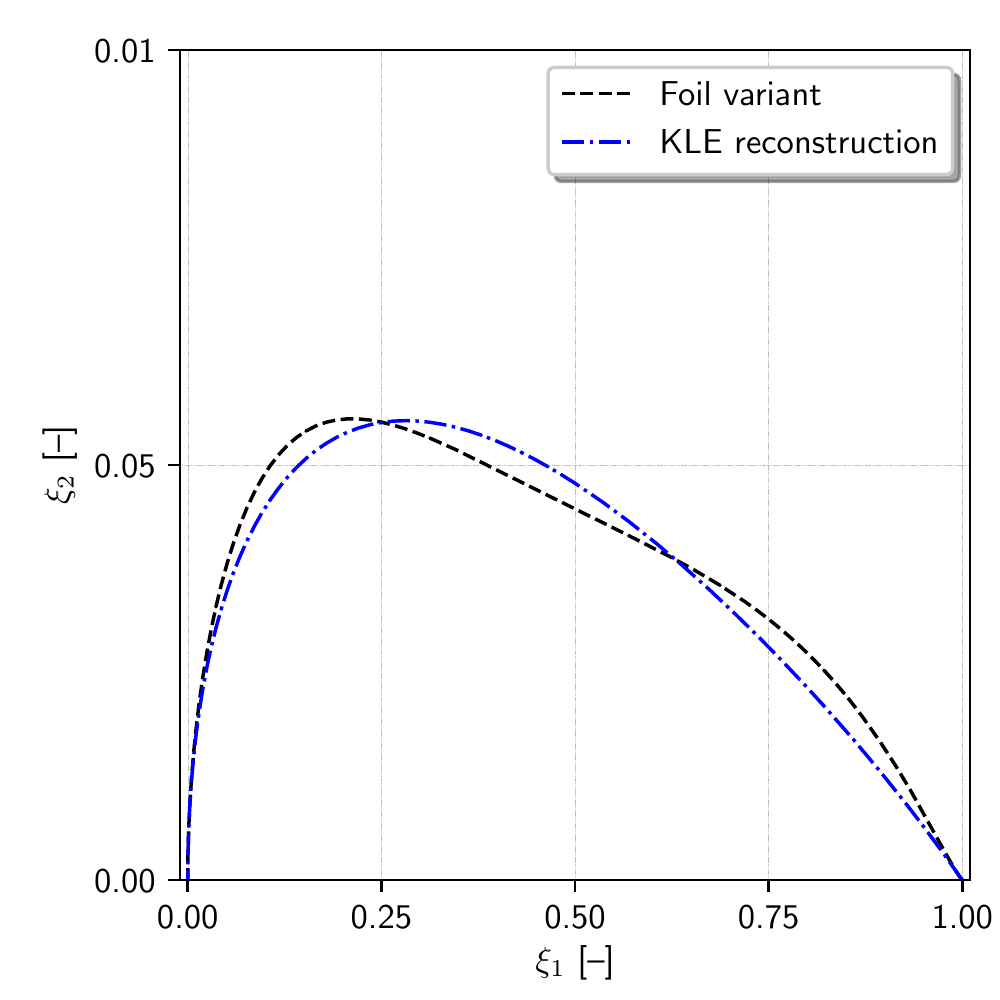}
\includegraphics[width=0.32\textwidth]{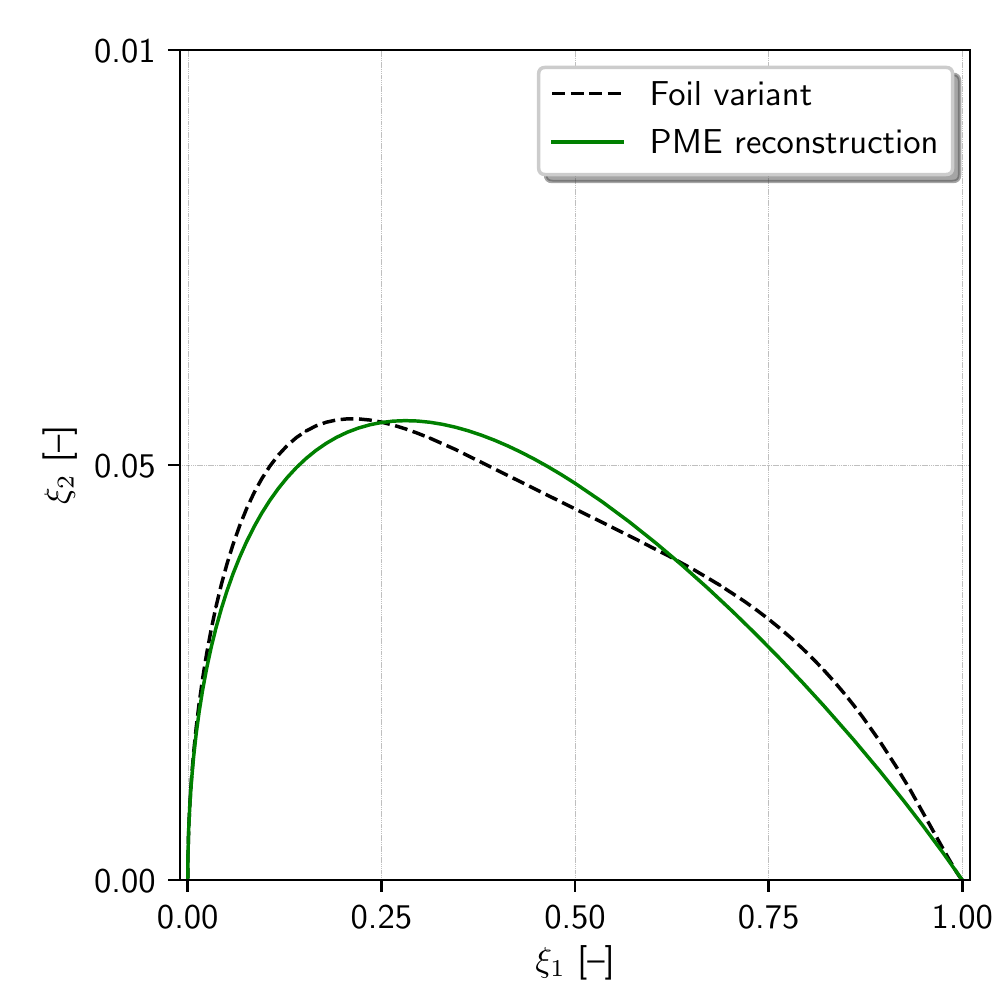}\\
\includegraphics[width=0.85\textwidth]{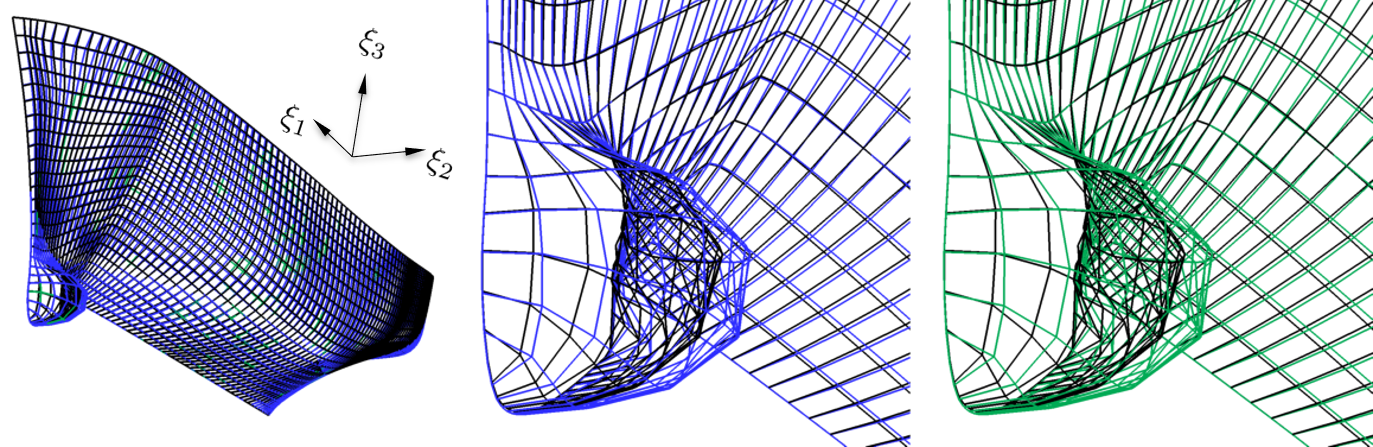}
\caption{Example of reconstruction of a geometry of the data set $\mathbf{D}$ \blue{for NACA 0012 (top) and DTMB 5415 (bottom) problems: shapes superposition (left) and details of differences (center and right)}}\label{fig:rec}
\end{figure}

Figure \ref{fig:eigsum} \change{(center)} shows the variance resolved by KLE and PME and the corresponding normalized mean squared error (NMSE, \change{see Fig. \ref{fig:eigsum} right}) in reconstructing the matrix $\mathbf{D}$, where the desired level of variance to be retained was set to 95\%. Specifically, Fig. \ref{fig:eigsum} (center) shows how KLE and PME cumulative sums of the eigenvalues (as percentage of the total variance) perfectly coincide. The number of reduced design variables to retain at least the 95\% of the original geometric variance is equal to \blue{$N = 4$ and 7, for NACA 0012 and DTMB 5415, respectively,} achieving a dimensionality reduction close to 70\% \blue{for both problems}. Fig. \ref{fig:eigsum} (right) provides the corresponding \blue{NMSE = 2.5 and 4\% (for NACA 0012 and DTMB 5415, respectively,)}, showing again the consistency of the PME approach with KLE. The corresponding eigenfunctions $\bfphi_k$, that can be used as shape modification basis by KLE, are shown in \blue{Figs. \ref{fig:eigg_foil} and} \ref{fig:eigg} (left). Note that the PME method provides exactly the same eigenvector components for the shape representation, but in addition provides also the basis functions $\bfv_k$ that embed the original design variables (see \blue{Figs. \ref{fig:eigg_foil} and} \ref{fig:eigg}, right), allowing to use original shape modification method/parameterization without the need of using the geometrical eigenfunctions $\bfphi_k$. \change{More in details, \blue{Figs. \ref{fig:eigg_foil} and} \ref{fig:eigg} (right) provides the basis $\mathbf{v}_k$ of PME eigenvectors used hereafter for the reconstruction of the original design parameterization as per Eq. \ref{eq:recu}.}
%

Once more, to show the consistency of PME with KLE, Fig. \ref{fig:corr} shows the correlation of the reconstruction error (normalized squared error, NSE) for each MC sample in $\mathbf{D}$, along with its PDF. It is evident how the two approaches perfectly match. Furthermore, it may be noted how \blue{(for both problems)} the mode of the NSE is lower than the mean value (NMSE), meaning that there are a large number of reconstructions below the NMSE threshold (complement of the variance retained). An example of geometry reconstruction is shown in Fig. \ref{fig:rec}. Specifically, Fig. \ref{fig:rec} \blue{(left)} shows one of the 1000 MC items in $\mathbf{D}$ obtained with \blue{the original design parameterization (Bezier on top and FFD on bottom) and the corresponding reconstruction via} KLE and PME. Minor differences can be seen comparing the \blue{original parameterization} (black) with KLE (blue) and PME (green) and details are provided in Fig. \ref{fig:rec} (centre and right). KLE in blue and PME in green show the same reconstruction error due to the truncation of the eigenvector expansion up to $N$ (see Eqs. \ref{eq:delta} and \ref{eq:recu}, respectively). 

\begin{figure}[!t]
\centering
\includegraphics[width=1\textwidth]{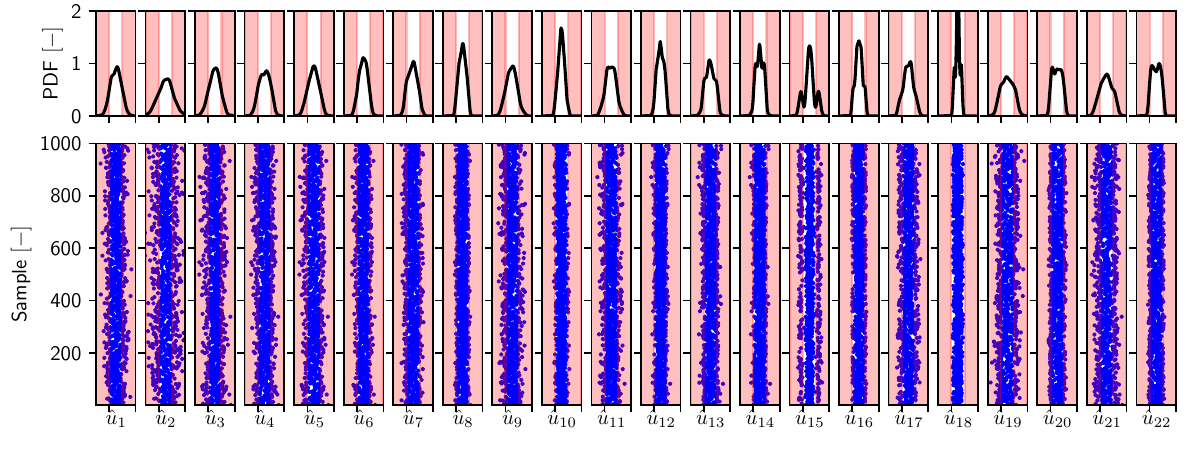}
\caption{Reconstruction of the original design variable via PME for 1000 random samples }\label{fig:bounds}
\end{figure}

Finally, \blue{considering only the DTMB 5415 problem (for the sake of simplicity) and} selecting 1000 random realizations of the reduced design variable vector $\bfx$, Fig. \ref{fig:bounds} shows the corresponding reconstruction of the original design variables $\widehat{\bfu}$: on bottom, each $j$-th sample is shown with a blue dot, for each component of $\widehat{\bfu}$ (columns); on top, the corresponding PDF. It can be noted, that embedding the design variables allows to fill the whole original design space, but (at the same time) provides many design variable vectors outsize the domain bounds \change{(shown as the light-red bands)}, thus producing possible unfeasibilities. This characteristics is common to both KLE and PME. Nevertheless, PME allows to penalize those design \change{that} fall outside the original design bounds.
\begin{figure}[!b]
\centering
\includegraphics[width=0.45\textwidth]{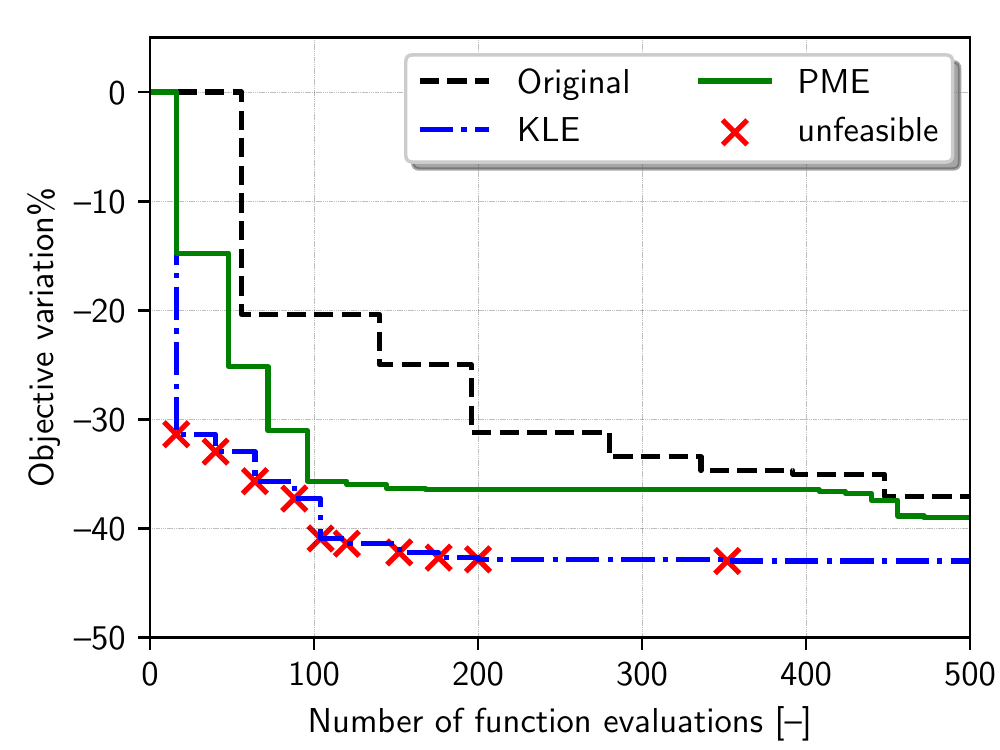}
\includegraphics[width=0.45\textwidth]{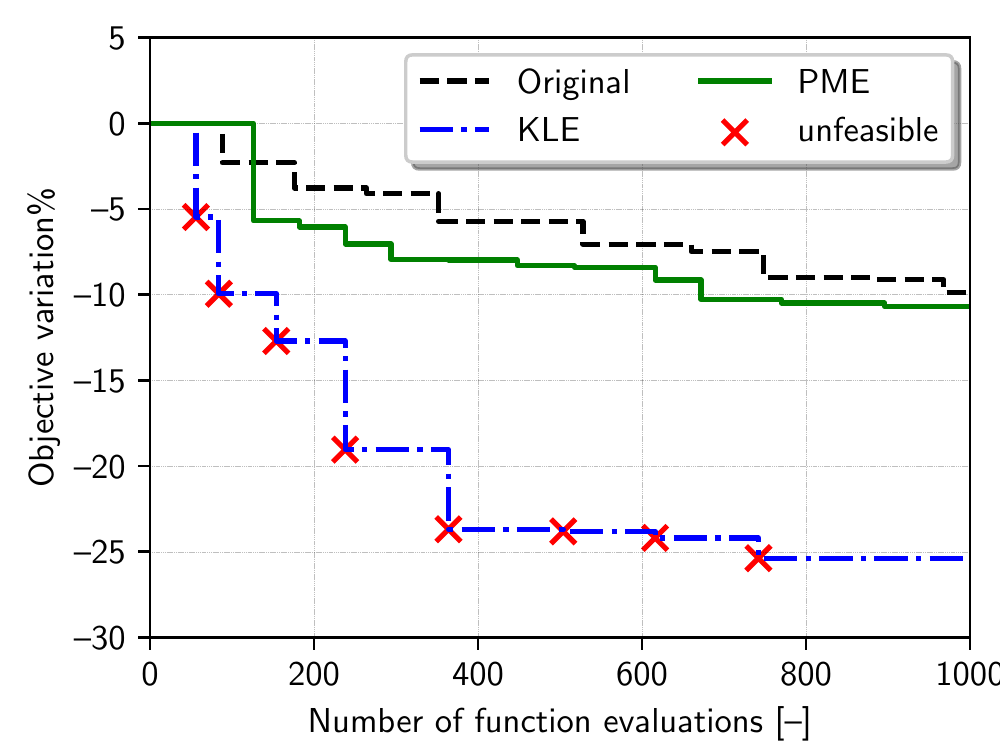}
\caption{Comparison of the optimization convergence \blue{for NACA 0012 (left) and DTMB 5415 (right) problems, }within original, KLE, and PME design spaces}\label{fig:opt}
\end{figure}

\subsection{Optimization results}
\blue{The solutions of problems in Eqs. \ref{eq:SDDOprob_foil} and \ref{eq:SDDOprob} are obtained} with a memetic version (hybrid global/local) \cite{serani2015-RASIEC} of the deterministic particle swarm optimization algorithm \cite{serani2016-ASC}. A total computational budget of \blue{500 and 1000 function evaluations is used for NACA 0012 and DTMB 5415 problems, respectively}. The optimization are conducted in the original design spaces, as well as using KLE and PME reduced design spaces.   

\change{For PME design space only, those designs that fall outside the original design space are penalized by a linear penalty function, proportional to distance between the design point and the design variable bounds. Specifically, the penalty is directly added to the objective function value $f(\bfu)$, used by the optimization algorithm, as follows
\begin{equation}
f(\bfu) = f(\bfu) + c\sum_{j=1}^M \max\left(u_j^l-u_j, u_j-u_j^u\right) \qquad \mathrm{if} \,\, (u_j<u_j^l \,\, \mathrm{or} \,\, u_j>u_j^u)   
\end{equation}
with $c=1000$.
}

\begin{figure}[!t]
\centering
\includegraphics[width=1\textwidth]{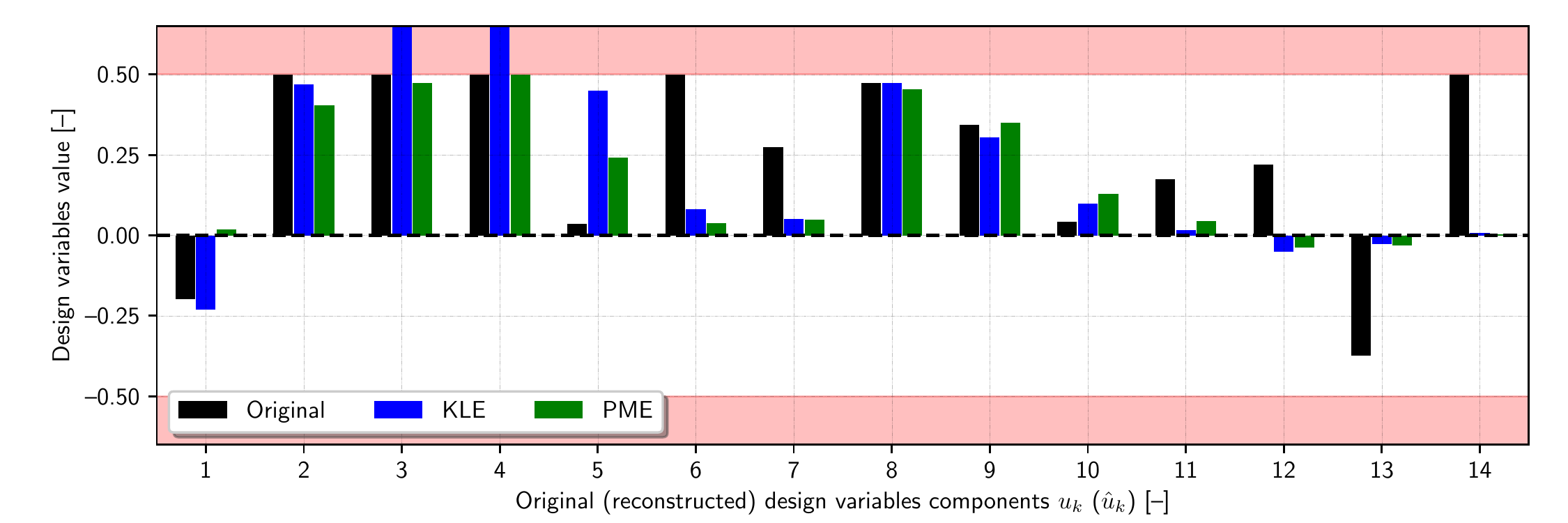}\\
\includegraphics[width=1\textwidth]{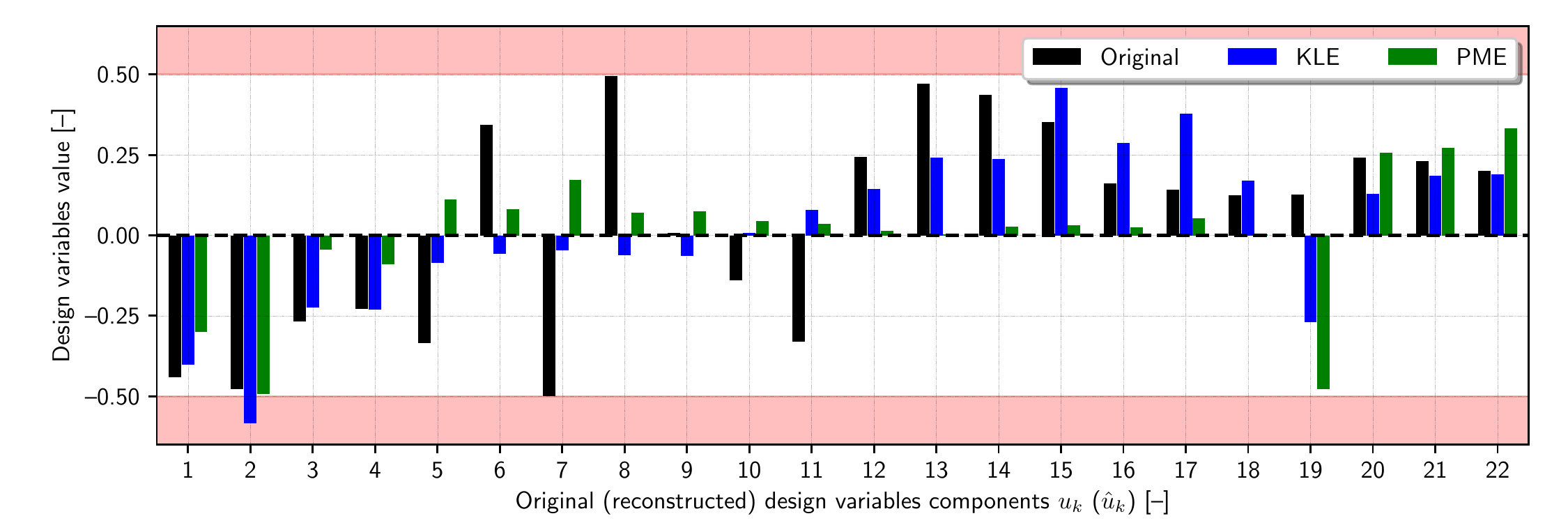}
\caption{Comparison of the optimal design variables \blue{for NACA 0012 (top) and DTMB 5415 (bottom) problems,} using the original and reduced design spaces}\label{fig:opt_var}
\end{figure}

Figure \ref{fig:opt} shows the optimization convergence \blue{for NACA 0012 (left) and DTMB 5415 (right)}. The optimization performed in the reduced design spaces outperforms the optimization procedure in the original design space, \blue{for both SDDO problems}. It may be noted how the optimization in the KLE design space (see blue dash-dotted line in Fig. \ref{fig:opt}) provides the best objective improvement, nevertheless all the optimal designs found during the optimization-algorithm evolution \blue{(for both SDDO problems)} fall outside the original design bounds (see red cross Fig. \ref{fig:opt}). This means that the KLE domain extends beyond the original \blue{parameterization} bounds, thus making the comparison unfair at least for the current applications. PME, on the one hand, spans the same KLE space, while, on the other hand, allows for easy penalization of unfeasible designs, effectively confining the optimization process within the original design bounds. 
\change{Comparing the optimization using the original and PME spaces, it can be noted how the major drops in the objective functions (about \blue{25, 30, and 35\% for NACA 0012 and} 5, 7.5, and 10\% \blue{for DTMB 5415}) are achieved around \blue{50, 75, and 100 (for NACA 0012) and} 150, 300, and 650 \blue{(for DTMB 5415)} function evaluations using PME, while similar drops are achieved in the original design space only after about \blue{150, 200, and 450 (for NACA 0012) and} 350, 650, and 1000 \blue{(for DTMB 5415)} function evaluations. Consequently, it can be concluded that PME, at least for the current applications, allows a better exploration and exploitation of the design space by the optimization algorithm, compared to the original design space, thus providing a faster convergence towards a better optimum, \blue{with an average (potential) computational saving of 60\%}.} 
\begin{figure}[!t]
\centering
\includegraphics[width=0.325\textwidth]{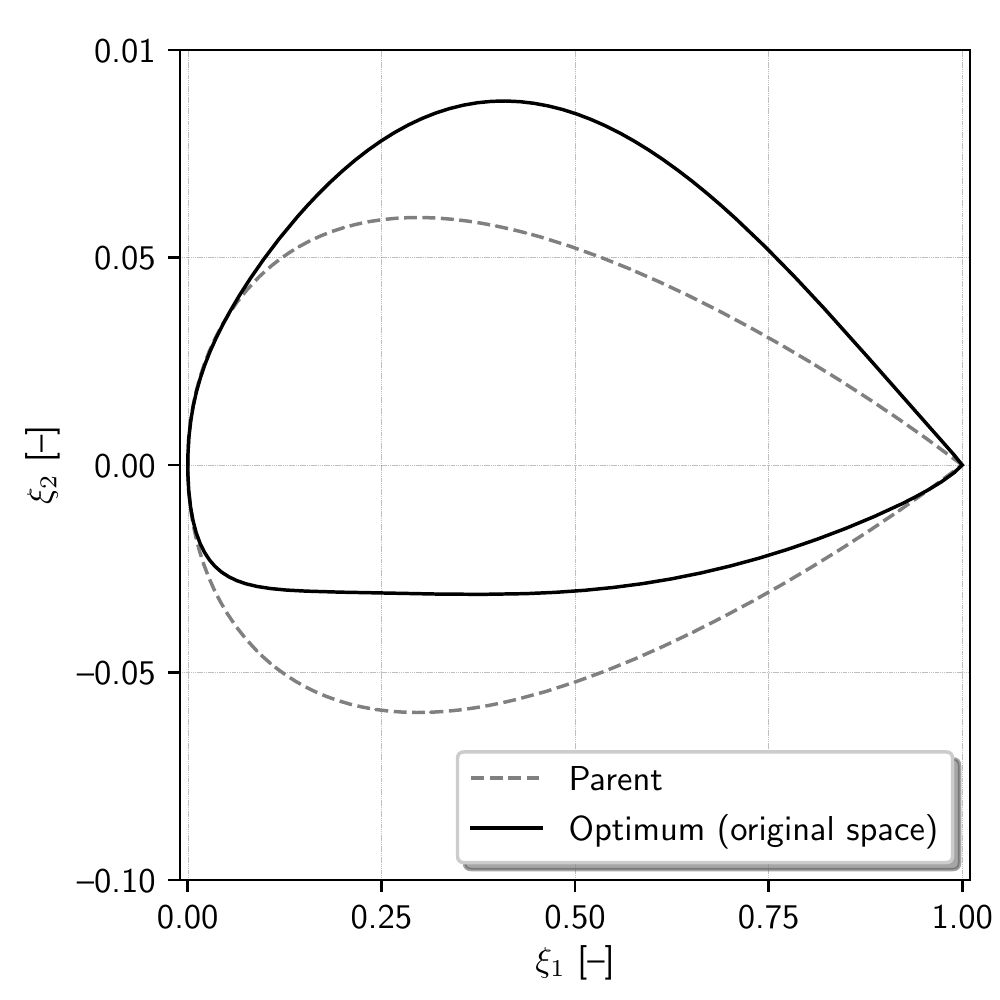}
\includegraphics[width=0.325\textwidth]{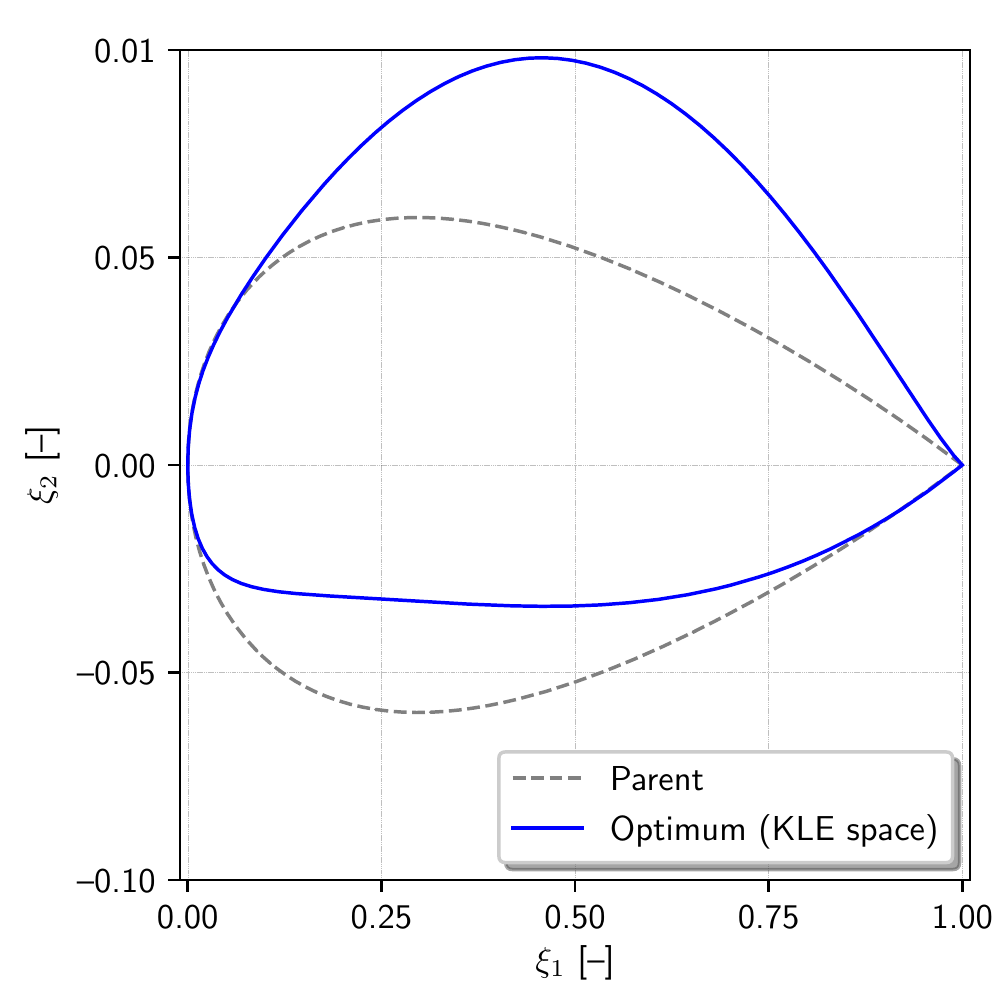}
\includegraphics[width=0.325\textwidth]{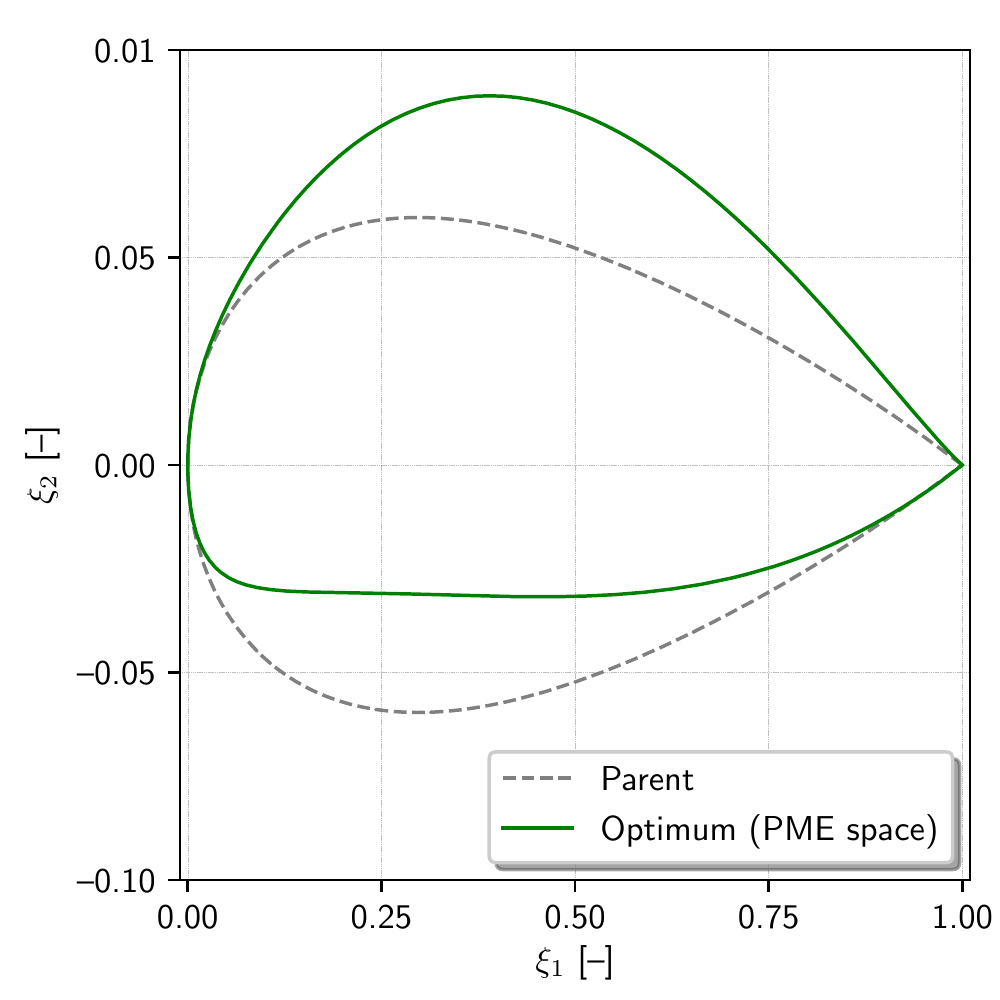}\\
\includegraphics[width=0.325\textwidth]{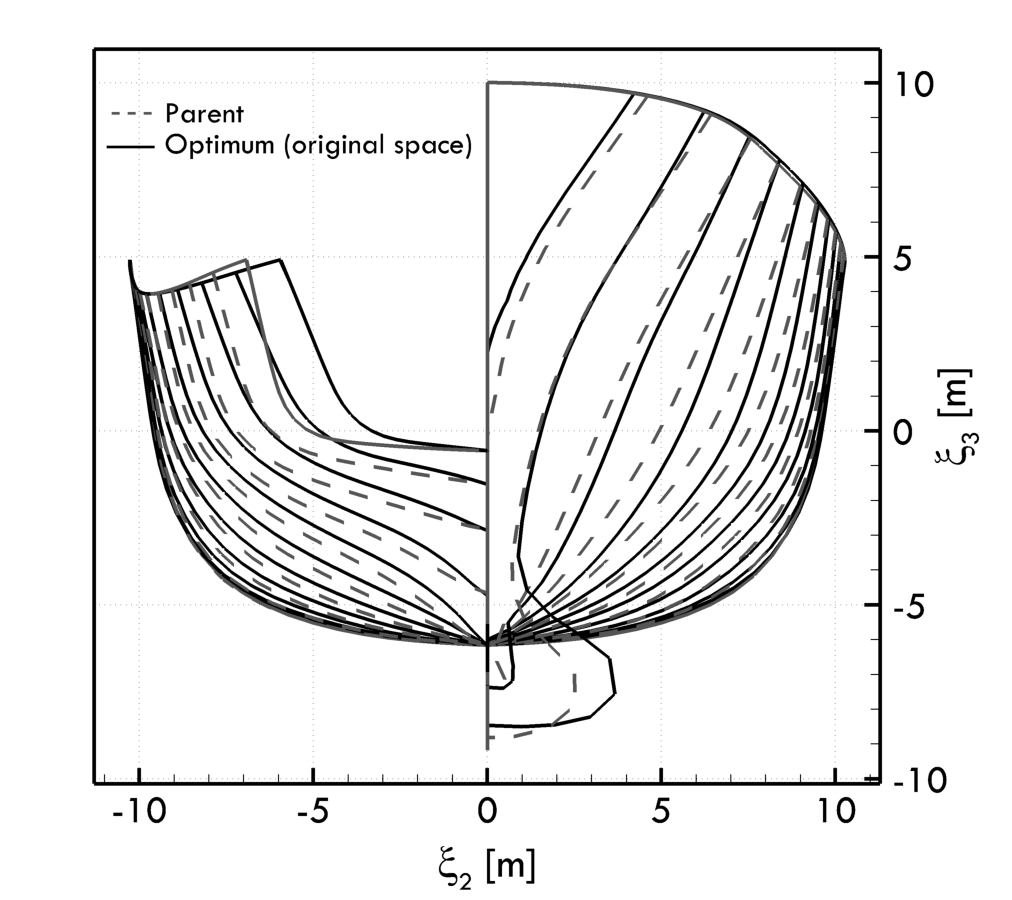}
\includegraphics[width=0.325\textwidth]{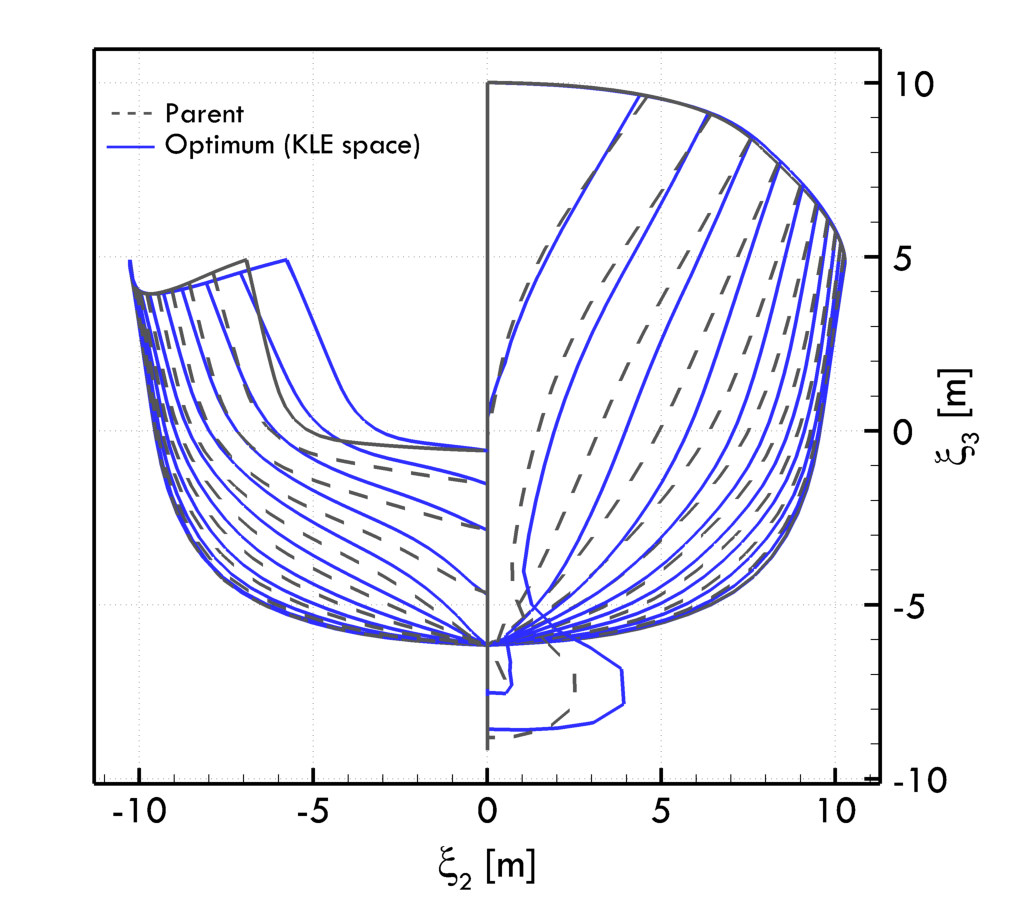}
\includegraphics[width=0.325\textwidth]{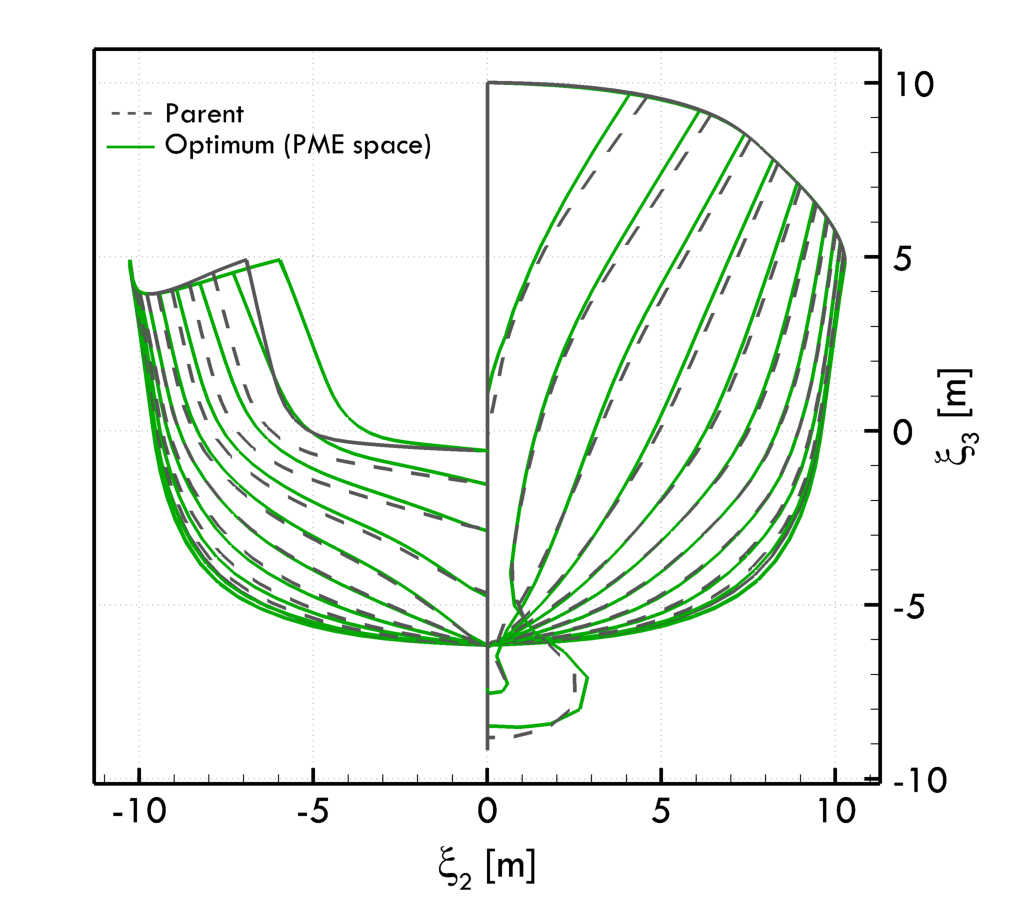}
\caption{\blue{Shape} comparison for the optimal designs \blue{of the NACA 0012 (top) and DTMB 5415 (bottom) problems}: from left to right original parameterization, KLE, and PME}\label{fig:opt_stations}
\end{figure}

\change{
\blue{Finally, Fig. \ref{fig:opt_var} provides a comparison of the optimal design variables for NACA 0012 (top) and DTMB 5415 (bottom) problems, respectively, giving the original parameterizations and their reconstruction from PME of the optimal design variants}. It can be seen that KLE violates \blue{the upper bound of third and fourth design variables for NACA 0012 and} the lower bound of the second design variable only \blue{for the DTMB 5415 problem}. All the optima have significantly different design variables, meaning that the optimization in the reduced-design spaces have allowed a better global exploration of the original design space, reaching better optima, located in different region of the \blue{original parameterizations} (note that only PME provides a feasible optimum). These differences reflects in the corresponding geometries, as shown in Fig. \ref{fig:opt_stations}, \blue{specially for the DTMB 5415 hull stations (iso-lines orthogonal to the hull longitudinal axis, see Fig. \ref{fig:opt_stations} bottom)}.}

\section{Conclusions and future work}\label{sec:V}
A methodology to address the curse of dimensionality in shape optimization has been presented. The parametric model embedding method has demonstrated its ability to provide a reduced dimensionality representation of the original design parameterization, using the original design parameters, allowing also for a faster optimization convergence towards a better optimum compared with the original design space. The method extends the original design-space dimensionality reduction procedure based on the KLE of the shape modification vector, presented in \cite{diez2015-CMAME}. The data matrix formed by realization of the shape modification vector is augmented with the associated design variables. To the latter is imposed a null weight. This provides to PME the capability of achieving the same results of KLE in terms of design-space assessment and dimensionality reduction. In addition, PME provides a basis to embed directly the original design variables. This last point has a notable industrial implication, because it overcomes one of the limitations associated to the original KLE formulation, i.e., the difficulty in returning to the original parameterization. 

PME is demonstrated for \blue{Bezier curves and} FFD-based design spaces of \blue{14 and} 22 design variables, \blue{respectively}, providing the pathway for a possible direct integration in CAD software and simulation-driven design optimization toolchain, avoiding the necessity to implement a new shape modification method based on the geometric component of the KLE eigenvectors and also easily enabling a search for the optimum within the original design bounds. \blue{Specifically, for the current applications, the PME has provided about 70\% design-space dimensionality reduction, achieving a 60\% computational budget saving, on average, to get at least the same objective function improvement obtained with the original design parameterization.}

Despite the evident benefits that PME provides to shape optimization, it is still relying on the linearity approximation of KLE and could not be as much efficient when strong nonlinearities are present in the shape parameterization or when the design-space dimensionality reduction procedure is extended to its physics-informed versions \cite{diez2016-SNH,serani2017-AIAA}, where nonlinearities can arise from the physical quantities. For these reasons, future work will investigate the possible extension of PME to nonlinear dimensionality reduction methods \cite{dagostino2017-MOD}, including studies on physics-informed formulations \cite{serani2020-JSR, khan2022geometric}.
\appendix
\change{
\section{Relationship between KLE and PME eigensolutions}\label{sec:appA}
Consider, for the sake of simplicity, $\mathbf{G}=\mathbf{I}$. Given the eigenproblems in Eqs. \ref{eq:pca}, its eigenvectors $\mathbf{Z}$ are orthonormal such that
\begin{equation}
\mathbf{ZWZ}^\mathsf{T}=\mathbf{I}    
\end{equation}
and the corresponding eigenvalues are those that nullify the determinant of the weighted autocovariance matrix of KLE as follows
\begin{equation}
\mathrm{det}\left(\mathbf{AW}-\lambda\mathbf{I}\right)=0    
\end{equation}
%
}

\change{
For the PME, starting from the definition 
\begin{equation}
\widetilde{\mathbf{A}}=\frac{1}{S}\mathbf{P}\mathbf{P}^\mathsf{T}=
\frac{1}{S}
\left[
\begin{array}{c}
\mathbf{D}\\
\mathbf{U}
\end{array}
\right]
\left[
\begin{array}{cc}
\mathbf{D}^\mathsf{T}&\mathbf{U}^\mathsf{T}
\end{array}
\right]
=
\left[
\begin{array}{cc}
\dfrac{1}{S}\mathbf{D}\mathbf{D}^\mathsf{T} & \dfrac{1}{S}\mathbf{D}\mathbf{U}^\mathsf{T}\\
\\
\dfrac{1}{S}\mathbf{U}\mathbf{D}^\mathsf{T} & \dfrac{1}{S}\mathbf{U}\mathbf{U}^\mathsf{T}\\
\end{array}
\right]
=
\left[
\begin{array}{cc}
\mathbf{A} & \dfrac{1}{S}\mathbf{D}\mathbf{U}^\mathsf{T}\\
\\
\dfrac{1}{S}\mathbf{U}\mathbf{D}^\mathsf{T} & \dfrac{1}{S}\mathbf{U}\mathbf{U}^\mathsf{T}\\
\end{array}
\right]
\end{equation}
and multiplying $\widetilde{\mathbf{A}}$ by the weight matrix $\widetilde{\mathbf{W}}$ in Eq. \ref{eq:pme_weight} leads to
\begin{equation}
\widetilde{\mathbf{A}}\widetilde{\mathbf{W}}=
\left[
\begin{array}{cc}
\mathbf{A} & \dfrac{1}{S}\mathbf{D}\mathbf{U}^\mathsf{T}\\
\\
\dfrac{1}{S}\mathbf{U}\mathbf{D}^\mathsf{T} & \dfrac{1}{S}\mathbf{U}\mathbf{U}^\mathsf{T}\\
\end{array}
\right]
\left[
\begin{array}{cc}
\mathbf{W} & \mathbf{0}\\
\mathbf{0} & \mathbf{0}\\
\end{array}
\right]
=
\left[
\begin{array}{cc}
\mathbf{AW} & \mathbf{0}\\
\mathbf{CW} & \mathbf{0}\\
\end{array}
\right]
\qquad
\mathrm{with}
\qquad
\mathbf{C}= \dfrac{1}{S}\mathbf{U}\mathbf{D}^\mathsf{T}
\end{equation}
The eigenvalues of the problems in Eq. \ref{eq:pme_pca} are given by
\begin{equation}
\mathrm{det}\left(\widetilde{\mathbf{A}}\widetilde{\mathbf{W}}-\lambda\mathbf{I}\right)=0    
\end{equation}
Using the Laplace expansion for the evaluation of the determinant of a generic square matrix $\mathbf{B}$
\begin{equation}
\mathrm{det}(\mathbf{B})=\sum_{j=1}^n (-1)^{i+j}\mathbf{B}_{i,j}\mathbf{M}_{i,j}    
\end{equation}
where $\mathbf{B}_{i,j}$ is the entry of the $i$-th row and $j$-th column of $\mathbf{B}$, and $\mathbf{M}_{i,j}$ is the determinant of the submatrix obtained by removing the $i$-th row and the $j$-th column of $\mathbf{B}$, and exploiting the null columns of $\widetilde{\mathbf{A}}\widetilde{\mathbf{W}}$ yields to
\begin{equation}
\mathrm{det}\left(\widetilde{\mathbf{A}}\widetilde{\mathbf{W}}-\lambda\mathbf{I}\right)=\left(\prod_{i=1}^M \lambda_i\right) \,\, \mathrm{det}\left({\mathbf{A}}{\mathbf{W}}-\lambda\mathbf{I}\right) 
\end{equation}
consequently the PME has many null eigenvalues as the number of design variables $M$ and the remaining are equivalent to those of KLE, such that
\begin{equation}\label{eq:aeig}
\widetilde{\boldsymbol{\Lambda}}=
\left[
\begin{array}{cc}
\boldsymbol{\Lambda} & \mathbf{0}\\
\mathbf{0} & \mathbf{0}\\
\end{array}
\right]
\end{equation}
Given Eq. \ref{eq:aeig}, the eigenproblem Eq. \ref{eq:pme_pca} can be written as follows
\begin{equation}\label{eq:eig_dem}
\left[
\begin{array}{cc}
\mathbf{AW} & \mathbf{0}\\
\mathbf{CW} & \mathbf{0}\\
\end{array}
\right]
\widetilde{\mathbf{Z}}
=
\widetilde{\mathbf{Z}}
\left[
\begin{array}{cc}
\boldsymbol{\Lambda} & \mathbf{0}\\
\mathbf{0} & \mathbf{0}\\
\end{array}
\right]
\end{equation}
Defining the eigenvectors $\widetilde{\mathbf{Z}}$ as the following generic block matrix
\begin{equation}
\widetilde{\mathbf{Z}}
=
\left[
\begin{array}{cc}
\mathbf{Q} & \mathbf{E}\\
\mathbf{V} & \mathbf{H}\\
\end{array}
\right]
\end{equation}
Eq. \ref{eq:eig_dem} becomes
\begin{equation}\label{eq:eig_dem}
\left[
\begin{array}{cc}
\mathbf{AWQ} & \mathbf{0}\\
\mathbf{CWQ} & \mathbf{0}\\
\end{array}
\right]
=
\left[
\begin{array}{cc}
\mathbf{Q}\boldsymbol{\Lambda} & \mathbf{0}\\
\mathbf{V}\boldsymbol{\Lambda} & \mathbf{0}\\
\end{array}
\right]
\end{equation}
yielding to
\begin{equation}\label{eq:eigeq}
\mathbf{AWQ}=\mathbf{Q}\boldsymbol{\Lambda} \implies \mathbf{Q} = \mathbf{Z}\\
\end{equation}
\begin{equation}
\mathbf{AWE}=\mathbf{0} \implies \mathbf{E} = \mathbf{0}\\
\end{equation}
\begin{equation}
\mathbf{CWQ}=\mathbf{V}\boldsymbol{\Lambda} \implies \mathbf{V} = \mathbf{CWZ}\boldsymbol{\Lambda}^{-1}
\end{equation}
consequently, also KLE and PME geometrical eigenvectors components are equivalent (see Eq. \ref{eq:eigeq}).
}
\section*{Acknowledgments}
The work is conducted in collaboration with the NATO task group AVT-331 on “Goal-driven, multi-fidelity approaches for military vehicle system-level design”. The authors are grateful to the US Office of Naval Research Global for its support through grants N62909-11-1-7011 and N62909-21-1-2042.

\bibliographystyle{abbrv}  
\bibliography{biblio}  

\end{document}